\documentclass{article}

\usepackage{amsmath}
\usepackage{amssymb}
\usepackage{theorem}
\usepackage{xypic}

\renewcommand{\tilde}{\widetilde}

\newcommand{\toto}{\rightrightarrows}

\newcommand{\MM}{{\mathfrak M}}

\newcommand{\Mm}{{\mathfrak m}}

\newcommand{\zz}{{\mathbb Z}}

\newcommand{\cc}{{\mathbb C}}

\newcommand{\rr}{{\mathbb R}}

\newcommand{\tT}{{\cal T}}

\renewcommand{\O}{{\cal O}}

\newcommand{\cC}{{\cal C}}
\newcommand{\eE}{{\cal E}}

\newcommand{\fF}{{\cal F}}

\newcommand{\hH}{{\cal H}}

\newcommand{\cl}{\mathop{\rm cl}\nolimits}

\newcommand{\tr}{\mathop{\rm tr}\nolimits}

\newcommand{\res}{\mathop{\rm res}\nolimits}

\newcommand{\sheafhom}{\mathop{\rm {\mit{ \hH\! om}}}\nolimits}

\newcommand{\sheafext}{\mathop{\rm {\mit{ \eE\! xt}}}\nolimits}
\newcommand{\sheaftor}{\mathop{\rm {\mit{ \tT\! or}}}\nolimits}

\newcommand{\comp}{\mathbin{{\scriptstyle\circ}}}
\newcommand{\ol}{\overline}

\newcommand{\doublearrowstack}[2]%
                      {{{{\scriptstyle#1}\atop{\textstyle\longrightarrow}}\atop{{\textstyle\longrightarrow}\atop{\scriptstyle#2}}}}
\newcommand{\rightleftarrowstack}[2]%
                      {{{{\scriptstyle#1}\atop{\textstyle\longrightarrow}}\atop{{\textstyle\longleftarrow}\atop{\scriptstyle#2}}}}
\newcommand{\leftrightarrowstack}[2]%
                      {{{{\scriptstyle#1}\atop{\textstyle\longleftarrow}}\atop{{\textstyle\longrightarrow}\atop{\scriptstyle#2}}}}

\newtheorem{thm}{Theorem}[section]
\newtheorem{cor}[thm]{Corollary}
\newtheorem{lem}[thm]{Lemma}
\newtheorem{prop}[thm]{Proposition}

\theorembodyfont{\upshape}
\newtheorem{defn}[thm]{Definition}

\newtheorem{rem}[thm]{Remark}
\newtheorem{rmk}[thm]{Remark}

\newtheorem{ex}[thm]{Example}

\newenvironment{pf}{\begin{trivlist}\item[]{\sc Proof.}}%
            {\nolinebreak $\Box$ \end{trivlist}}

\newenvironment{proof}{\begin{trivlist}\item[]{\sc Proof.}}%
            {\nolinebreak $\Box$ \end{trivlist}}

\newcommand{\qcoh}{\text{\rm Qcoh-}}
\newcommand{\vir}{{\rm vir}}

\newcommand{\A}{\mathbb A}

\newcommand{\Gm}{{{\mathbb G}_m}}

\newcommand{\cpx}{\mathbb C}

\newcommand{\eps}{\varepsilon}

\newcommand{\del}{\partial}
\newcommand{\longto}{\longrightarrow}
\DeclareMathOperator\id{id}
\newcommand{\ltensor}{\mathbin{\displaystyle\stackrel{L}{\otimes}}}

\renewcommand{\O}{\mathcal O}

\newcommand{\cok}{\mathop{\rm cok}\nolimits}

\newcommand{\rk}{\mathop{\rm rk}\nolimits}

\newcommand{\spec}{\mathop{\rm Spec}\nolimits}
\newcommand{\Hilb}{\mathop{\rm Hilb}\nolimits}

\newcommand{\Hom}{\mathop{\rm Hom}\nolimits}

\newcommand{\noprint}[1]{}

\def\Label#1{\label{#1}{\tt [#1]}\phantom{h}}

\newcommand{\unsure}[1]{{\footnotesize #1}}

\def\Label{\label}
\def\unsure{\noprint}

\author{K. Behrend and B. Fantechi}
\title{Symmetric Obstruction Theories and\\
Hilbert Schemes of Points on Threefolds}

\begin{document}
\sloppy

\maketitle

\abstract{Recall that in an earlier paper by one of the authors
  Donaldson-Thomas type invariants were expressed as certain weighted Euler
  characteristics of the moduli space. The Euler characteristic is
  weighted by a certain canonical $\zz$-valued constructible function
  on the moduli space. This constructible function associates to any
  point of the moduli space a certain invariant of the singularity of
  the space at the point.

In the present paper, we evaluate this invariant for the case of a
singularity which is an isolated point of a $\cc^\ast$-action and
which admits a symmetric obstruction theory compatible with the
$\cc^\ast$-action. The answer is $(-1)^d$, where $d$ is the dimension
of the Zariski tangent space. 

We use this result to prove that for any threefold, proper or not, the
weighted Euler characteristic of the Hilbert scheme of $n$ points on
the threefold is, up to sign, equal to the usual Euler
characteristic. For the case of a projective Calabi-Yau threefold,
we deduce that the Donaldson-Thomas invariant of the Hilbert scheme
of $n$ points is, up to sign, equal to the Euler characteristic.  This
proves a conjecture of Maulik-Nekrasov-Okounkov-Pandharipande.}

\eject
\tableofcontents

\eject
\section*{Introduction}
\addcontentsline{toc}{section}{Introduction}

\subsubsection{Symmetric obstruction theories}

The first purpose of this paper is to introduce {\em symmetric
  obstruction theories}. In a nutshell, these are obstruction theories
for which the space of infinitesimal deformations is the dual of the
space of infinitesimal obstructions. 

As an example of an obstruction theory, consider the intersection of
two smooth varieties $V$, $W$ inside an ambient smooth variety
$M$. The intersection $X$ carries an obstruction theory. This is the
2-term complex of vector bundles 
$$\xymatrix@1{E\ar@{}[r]|-=&[\,\Omega_M|_X\rrto^-{\res_V-\res_W}&&
\Omega_V|_X\oplus\Omega_W|_X\,]\,}$$ 
considered as an object of the derived category $D(X)$ of $X$, taking
up degrees $-1$ and $0$. We see that infinitesimal deformations of $X$ are
classified by $h^0(E^\vee)=T_X$, the sheaf of derivations on $X$.
Moreover, the obstructions to the smoothness of $X$ are contained in
$h^1(E^\vee)$, which is called the {\em obstruction sheaf}, notation
$ob=h^1(E^\vee)$. Note that $h^0(E^\vee)$ is intrinsic to $X$, but
$h^1(E^\vee)$ is not.  In fact, if $X$ is smooth, all obstructions
vanish, but $h^1(E^\vee)$ may be non-zero, although it is always a
vector bundle, in this case.

This obstruction theory $E$ is {\em symmetric}, if $M$ is a complex
symplectic manifold, i.e., hyperk\"ahler, and $V$, $W$ are Lagrangian
submanifolds. In fact, the symplectic form $\sigma$ induces a
homomorphism $T_X\to \Omega_M$, defined by $v\mapsto \sigma(v,-)$.
The fact that $V$ and $W$ are Lagrangian, i.e., equal to their own
orthogonal complements with respect to $\sigma$, implies that there is
an exact sequence
$$\xymatrix@1{
0\rto & T_X\rto & \Omega_M|_X\rto & \Omega_V|_X\oplus\Omega_W|_X \rto
& \Omega_X\rto &0\,.}$$
Hence, assuming for simplicity that $X$ is smooth and hence this is an
exact sequence of vector bundles, we see that
$ob=h^1(E^\vee)=\Omega_X$ and hence $T_X$ is, indeed, dual to $ob$. 

In more abstract terms, what makes an obstruction theory $E$ 
symmetric, is a non-degenerate symmetric bilinear form of degree 1
$$\beta:E\ltensor E\longrightarrow \O_X[1]\,.$$

If $M$ is an arbitrary smooth scheme, then $\Omega_M$ is a symplectic
manifold in a canonical way, and the graph of any closed 1-form
$\omega$ is a Lagrangian submanifold.  Thus the scheme theoretic zero
locus  $X=Z(\omega)$ of $\omega$ is an example of the above, the
second Lagrangian being the zero section. 

As a special case of this, we may consider the Jacobian locus
$X=Z(df)$ of a regular function on a smooth variety $M$. It is endowed
with a canonical symmetric obstruction theory. In Donaldson-Thomas
theory, where the moduli space is heuristically the critical locus of
the holomorphic Chern-Simons functional, there is a canonical
symmetric obstruction theory, see~\cite{RT}. 

Unfortunately, we are unable to prove that every symmetric obstruction
theory is locally given as the zero locus of a closed 1-form on a
smooth scheme, even though we see no reason why this should not be
true. 

The best we can do is to show that the most general local example of a
symmetric obstruction theory is the zero locus of an {\em almost
  closed }1-form on a smooth scheme. A form $\omega$ is almost closed
if its differential $d\omega$ vanishes on the zero locus $Z(\omega)$.

For the applications we have in mind we also need equivariant versions
of all of the above, in the presence of a $\Gm$-action.

\subsubsection{Weighted Euler characteristics and $\Gm$-actions}

In \cite{Beh} a new (as far as we can tell) invariant of singularities
  was introduced.
For a singularity $(X,P)$ the notation was
$$\nu_X(P)\,.$$
The function $\nu_X$ is a constructible $\zz$-valued function on any
Deligne-Mumford stack $X$. In \cite{Beh}, the following facts were
proved about $\nu_X$:

\smallskip\noindent
$\bullet$ If $X$ is smooth at $P$, then $\nu_X(P)=(-1)^{\dim X}$. 

\smallskip\noindent
$\bullet$ $\nu_X(P)\,\nu_Y(Q)=\nu_{X\times Y}(P,Q)$. 

\smallskip\noindent
$\bullet$ If $X=Z(df)$ is the singular locus of a regular function $f$
on a smooth variety $M$, then $$\nu_X(P)=(-1)^{\dim M}(1-\chi(F_P))\,,$$
where $F_P$ is the Milnor fibre of $f$ at $P$. 

\smallskip\noindent
$\bullet$ Let  $X$ be a projective scheme endowed with a symmetric obstruction
theory. The associated Donaldson-Thomas type invariant (or virtual
count) is the degree of the associated virtual fundamental class. 
In this case, $\nu_X(P)$ is the contribution of the point $P$ to the
Donaldson-Thomas type invariant, in the sense that
$$\#^\vir(X)=\chi(X,\nu_X)=\sum_{n\in\zz}n\,\chi(\{\nu_X=n\})\,.$$
We define the weighted Euler characteristic of $X$ to be
$$\tilde\chi(X)=\chi(X,\nu_X)\,.$$

The last property shows the importance of $\nu_X(P)$ for the
calculation of Donaldson-Thomas type invariants. 

In this paper we calculate the number $\nu_X(P)$ for certain kinds of
singularities. In fact, we will assume that $X$ admits a $\Gm$-action
and a symmetric obstruction theory, which are compatible with each
other.  Moreover, we assume $P$ to be an isolated fixed point for the
$\Gm$-action.  We prove that
\begin{equation}\Label{bigdeal}
\nu_X(P)=(-1)^{\dim T_X|_P}\,,
\end{equation}
in this case. 

We get results of two different flavors from this:

\smallskip\noindent
$\bullet$ If the scheme  $X$ admits a globally defined 
$\Gm$-action with isolated fixed points and around every fixed point
admits a symmetric obstruction theory compatible with the
$\Gm$-action we obtain
\begin{equation}\Label{useful}
\tilde\chi(X)=\sum_P (-1)^{\dim T_X|_P}\,,
\end{equation}
the sum extending over the fixed points of the $\Gm$-action.
This is because non-trivial $\Gm$-orbits
do not contribute, the Euler characteristic of $\Gm$ being zero, and
$\nu_X$ being constant on such orbits.

\smallskip\noindent
$\bullet$  If $X$ is projective, with globally defined 
$\Gm$-action and symmetric obstruction theory, these two structures
being compatible, we get
\begin{equation}\Label{amusing}
\#^\vir(X)=\tilde\chi(X)=\sum_P(-1)^{\dim T_X|_P}\,.
\end{equation}

\subsubsection{An example}

It may be worth pointing out how to prove (\ref{bigdeal}) in a special
case. Assume the multiplicative group $\Gm$ acts on affine $n$-space
$\A^n$ in a linear way with non-trivial weights
$r_1,\ldots,r_n\in\zz$, so that the origin $P$ is an isolated fixed
point.  Let $f$ be a regular function on $\A^n$, which is invariant
with respect to the $\Gm$-action. This means that $f(x_1,\ldots,x_n)$
is of degree zero, if we assign to $x_i$ the degree $r_i$. 
Let $X=Z(df)$ be the scheme-theoretic critical set of $f$. The scheme
$X$ inherits a $\Gm$-action.  It also carries a symmetric obstruction
theory which is compatible with the $\Gm$-action.  

Assume that $f\in (x_1,\ldots,x_n)^3$. This is not a serious
restriction. It ensures that 
$T_X|_P=T_{\A^n}|_P$ and hence that $\dim T_X|P=n$.  

Let $\epsilon\in\rr$, $\epsilon>0$ and $\eta\in\cc$, $\eta\not=0$ be
chosen such that the Milnor fiber of $f$ at the origin may be defined
as 
$$F_P=\{x\in \cc^n\mid\text{$f(x)=\eta$ and $|x|<\epsilon$}\}\,.$$
It is easy to check that $F_P$ is invariant under the $S^1$-action on
$\cc^n$ induced by our $\Gm$-action. Moreover, the $S^1$-action on
$F_P$ has no fixed points.  This implies immediately that
$\chi(F_P)=0$ and hence that $\nu_X(P)=(-1)^n$. 

Even though we consider this example $(Z(df),P)$ to be the prototype
of a singularity admitting compatible $\Gm$-actions and symmetric
obstruction theories, we cannot prove that every such singularity is
of the form $(Z(df),P)$.  We can only prove that a singularity with
compatible $\Gm$-action and symmetric obstruction theory looks like
$(Z(\omega),P)$, where $\omega$ is an  almost closed $\Gm$-invariant
1-form on $\A^n$, rather than the { exact }invariant 1-form
$df$. This is why the proof of (\ref{bigdeal}) is more involved, in
the general case.
Rather than using the Milnor fiber, we use the expression of
$\nu_X(P)$ as a linking number, Proposition~4.22 of~\cite{Beh}. 

\subsubsection{Lagrangian intersections}

One amusing application of (\ref{amusing}) is the following formula.
Assume $M$ is a complex symplectic manifold with a Hamiltonian
$\cc^\ast$-action, all of whose fixed points are isolated. Let $V$ and
$W$ be invariant Lagrangian submanifolds. Assume their intersection is
compact. Finally, assume that the Zariski tangent space of the
intersection at every fixed point is even-dimensional. Then we can
express the intersection number as an Euler characteristic:
$$\deg([V]\cap[W])=\chi(V\cap W)\,.$$

\subsubsection{Hilbert schemes}

Our result is a
powerful tool for computing weighted Euler characteristics.  It is a
replacement for the lacking additivity of $\tilde\chi$ over strata.

As an example of the utility of (\ref{bigdeal}), we will
show in this paper that 
\begin{equation}\Label{alsonotbad}
\tilde\chi(\Hilb^n Y)=(-1)^n\chi(\Hilb^n Y)\,,
\end{equation}
for every smooth scheme $Y$ of dimension~3.

In particular, if $Y$ is projective and Calabi-Yau (i.e., has a chosen
trivialization $\omega_Y=\O_Y$), we get that
$$\#^\vir(\Hilb^nY)=(-1)^n\chi(\Hilb^nY)\,,$$
where $\#^\vir$ is the virtual count \`a la
Donaldson-Thomas~\cite{RT}. This latter formula was conjectured by
Maulik-Nekrasov-Okounkov-Pandharipande \cite{MNOP}.
Using the McMahon function $M(t)=\prod_{n=1}^\infty (1-t^n)^{-n}$, we
can also express this result as
$$\sum_{n=0}^\infty\#^\vir(\Hilb^nY)\,t^n=M(-t)^{\chi(Y)}\,.$$

The strategy for proving (\ref{alsonotbad}) is as follows. We first
consider the open Calabi-Yau threefold $\A^3$.  We exploit a suitable
$\Gm$-action on $\A^3$ to prove (\ref{alsonotbad}) for $Y=\A^3$, using
Formula~(\ref{useful}). At this point, we can drop all Calabi-Yau
assumptions.  

Let $F_n$ be the {\em punctual }Hilbert scheme.  It
parameterizes subschemes of $\A^3$ of length $n$ which are entirely
supported at the origin. Let $\nu_n$ be the restriction of
$\nu_{\Hilb^n \A^3}$ to $F_n$.  Formula (\ref{alsonotbad}) for
$Y=\A^3$ is equivalent to
\begin{equation}\Label{punctual}
\chi(F_n,\nu_n)=(-1)^n\chi(F_n)\,.
\end{equation}

Finally, using more or less standard stratification arguments, 
we express $\tilde\chi(\Hilb^n Y)$ in terms of $\chi(F_n,\nu_n)$. 
This uses the fact that the punctual Hilbert scheme
of $Y$ at a point $P$ is isomorphic to $F_n$. 
Then (\ref{punctual}) implies (\ref{alsonotbad}).

\subsubsection{Conventions}

We will work over the field of complex numbers. All stacks will be of
Deligne-Mumford type.  All schemes and stacks will be of finite type
over $\cc$. Hence the derived category $D_{\text{qcoh}}(\O_X)$, of
complexes of sheaves of $\O_X$-modules with quasi-coherent cohomology
is equivalent to the derived category $D(\qcoh\O_X)$ of the category
of quasi-coherent $\O_X$-modules (see Proposition~3.7 in Expos\'e~II
of SGA6).  To fix ideas, we will denote by $D(X)$ the latter derived
category and call it the derived category of $X$.  We will often write
$E\otimes F$ instead of $E\ltensor F$, for objects $E,F$ of $D(X)$.

Let $X$ be a Deligne-Mumford stack.  We will write $L_X$ for the
cutoff at $-1$ of the cotangent complex of $X$.  Thus, if $U\to X$ is
\'etale and $U\to M$ a closed immersion into a smooth scheme $M$, we
have, {\em canonically},
$$L_X|_U=[I/I^2\to\Omega_M|_X]\,,$$
where $I$ is the ideal sheaf of $U$ in $M$ and we think of the
homomorphism $I/I^2\to \Omega_M|_X$ of coherent sheaves on $U$ as a
complex concentrated in the interval $[-1,0]$. We will also call $L_X$
the  cotangent complex of $X$, and hope the reader will forgive this
abuse of language. The cotangent complex $L_X$ is an object of $D(X)$.

We will often use homological notation for objects in the derived
category.  This means that $E_n=E^{-n}$, for a complex $\ldots\to
E^{i}\to E^{i+1}\to\ldots$ in $D(X)$.

For a complex of sheaves $E$, we denote the cohomology sheaves by
$h^i(E)$. 

Let us recall a few sign conventions: If
$E=[E_1\stackrel{\alpha}{\longrightarrow} E_0]$ is a complex 
concentrated in the interval $[-1,0]$, then
$E^\vee=[E_0^\vee\stackrel{-\alpha^\vee}{\longrightarrow} E_1^\vee]$
is a complex concentrated in the interval $[0,1]$. Thus the shifted
dual $E^\vee[1]$ is given by
$E^\vee[1]=[E_0^\vee\stackrel{\alpha^\vee}{\longrightarrow} E_1^\vee]$
and concentrated, again, in the interval $[-1,0]$. 

If $\theta:E\to F$ is a homomorphism of complexes concentrated in the
interval $[-1,0]$, such that
$\theta=(\theta_1,\theta_0)$, then the shifted
dual $\theta^\vee[1]:F^\vee[1]\to E^\vee[1]$ is given by
$\theta^\vee[1]=(\theta_0^\vee,\theta_1^\vee)$. 

Suppose $E=[E_1\stackrel{\alpha}{\longrightarrow}E_0]$ and
$F=[F_1\stackrel{\beta}{\longrightarrow} F_0]$ are complexes
concentrated in the interval $[-1,0]$ and  
$\theta:E\to F$ and $\eta:E\to F$ homomorphisms of complexes. Then  a
homotopy from $\eta$ to 
$\theta$ is a homomorphism $h:E_0\to F_1$ such that
$h\circ\alpha=\theta_1-\eta_1$ and $\beta\circ h=\theta_0-\eta_0$.

\subsubsection{Acknowledgments}

We would like to thank Jim Bryan and Lothar G\"ottsche for helpful
discussions.  We would especially like to thank Bumsig Kim and the Korea
Institute for Advanced Study for the warm hospitality during a stay in
Seoul, where the main results of this paper were proved.

\eject
\section{Symmetric Obstruction Theories}

\subsection{Non-degenerate symmetric bilinear forms}

\begin{defn}
Let $X$ be a scheme, or a Deligne-Mumford stack. Let $E\in
D^b_{coh}(\O_X)$ be a perfect 
complex.  A {\bf non-degenerate symmetric bilinear form of degree 1 }on $E$ is
a morphism 
$$\beta: E \ltensor E \longrightarrow \O_X[1]$$
in $D(X)$, which is 

(i) {\em symmetric}, which means that $$\beta(e\otimes
e')=(-1)^{\deg(e)\deg(e')}\beta(e'\otimes e)\,,$$ 

(ii) {\em non-degenerate}, which means that $\beta$ induces an isomorphism 
$$\theta:E\longrightarrow E^\vee[1]\,.$$
\end{defn}

\begin{rmk}
The isomorphism $\theta:E\to E^\vee[1]$ determines $\beta$ as the
composition
$$\xymatrix@1{E\otimes E\rto^-{\theta\otimes \id}&
E^\vee[1]\otimes E\rto^-{\tr[1]}& \O_X[1]\,.}$$
Symmetry of
$\beta$ is equivalent to the condition
$$\theta^\vee[1]=\theta\,.$$
Usually, we will find it more convenient to work with $\theta$, rather
than $\beta$. Thus we will think of a non-degenerate symmetric
bilinear form of degree 1 on
$E$ as an isomorphism $\theta:E\to E^\vee[1]$, satisfying
$\theta^\vee[1]=\theta$.
\end{rmk}

\begin{rmk}
Above, we have defined non-degenerate symmetric bilinear forms of
degree 1.  One 
can generalize the definition to any degree $n\in\zz$.  Only the case
$n=1$ will interest us in this paper.
\end{rmk}

\begin{ex}\Label{simple}
Let $F$ be a vector bundle on $X$ and let  $\alpha:F\to F^\vee$ a
symmetric bilinear form. Define the
complex $E=[ F\to F^\vee]$, by putting $F^\vee$ in degree $0$ and $F$ in degree
$-1$. Then $E^\vee[1]=E$. Define 
$\theta=(\theta_1,\theta_0)$ by 
$\theta_{1}=\id_{F}$ and $\theta_0=\id_{F^\vee}$: 
$$\xymatrix{
E\dto^\theta \ar@{}[r]|{=} & [F\rto^\alpha\dto_{1} & F^\vee\dto^1]\\
E^\vee[1] \ar@{}[r]|{=} & [F\rto^{\alpha} & F^\vee]}$$
Then $E$ is a  perfect
complex with perfect  amplitude contained in $[-1,0]$. Moreover,
$\theta$ is a
non-degenerate symmetric bilinear form on $E$. Note that $\theta$ is
an isomorphism, and hence the form it defines is non-degenerate,
whether or not $\alpha$ is non-degenerate. 
\end{ex}

\begin{ex}\Label{hess}
Let $f$ be a regular function on a smooth variety $M$.
The Hessian of $f$ defines a symmetric bilinear form on
$T_M|_X$, where $X=Z(df)$. So there is an induced symmetric
bilinear form
on the complex $[T_M|_X\to\Omega_M|_X]$.
\end{ex}

\begin{lem}\Label{jfdsk}
Let $E$ be a complex of vector bundles on $X$, concentrated in the
interval $[-1,0]$. Let $\theta:E\to 
E^\vee[1]$ be a homomorphism of complexes.  Assume that 
$\theta^\vee[1]=\theta$, as morphisms in the derived category. Then
Zariski-locally on the scheme $X$ (or \'etale 
locally on the stack $X$) we can represent the derived category
morphism given by $\theta$ as a homomorphism of
complexes  $(\theta_1,\theta_0)$:
$$\xymatrix{ E\dto^\theta\ar@{}[r]|{=} & [E_1\dto_{\theta_1}\rto^\alpha
  & E_0\dto^{\theta_0}]\\ 
E^\vee[1]\ar@{}[r]|{=} & [E_0^\vee \rto^{\alpha^\vee}& E_1^\vee]}$$
such that $\theta_1=\theta_0^\vee$.
\end{lem}
\begin{pf}
Let us use notation 
$\theta=(\psi_1,\psi_0)$. Then the equality of derived category
morphisms $\theta^\vee[1]=\theta$ implies
that, locally, $\theta^\vee[1]=(\psi_0^\vee, \psi_1^\vee)$ and
$\theta=(\psi_1, \psi_0)$ are homotopic.  So let $h:E_0\to
E_0^\vee$ be a homotopy: 
\begin{align*}
h\alpha & =\psi_1-\psi_0^\vee\\
\alpha^\vee h &=\psi_0-\psi_1^\vee\,.
\end{align*}
Now define 
\begin{align*}
\theta_0&={\textstyle\frac{1}{2}}(\psi_0+\psi_1^\vee)\\
\theta_1&={\textstyle\frac{1}{2}}(\psi_1+\psi_0^\vee)\,.
\end{align*}
One checks, using $h$,  that $(\theta_1,\theta_0)$ is a homomorphism
of complexes, and as 
such, homotopic to $(\psi_1,\psi_0)$. Thus $(\theta_1,\theta_0)$
represents the 
derived category morphism $\theta$, and has the required property.
\end{pf}

The next lemma shows that for amplitude 1 objects, every
non-degenerate symmetric bilinear form locally looks like the one given in
Example~\ref{simple}.  Again, locally means \'etale locally, but in
the scheme case Zariski locally.

\begin{lem}\Label{even.better}
Suppose that $A\in D^b_{coh}(\O_X)$ is of perfect amplitude contained
in $[-1,0]$, and that 
$\eta:A\to A^\vee[1]$ is an isomorphism satisfying $\eta^\vee[1]=\eta$.  
Then we can locally represent $A$ by a homomorphism of vector bundles
$\alpha:E\to E^\vee$ satisfying $\alpha^\vee=\alpha$ and the isomorphism
$\eta$ by the identity.
\end{lem}
\begin{pf}
Start by representing the derived category object $A$ by an actual complex
of vector bundles $\alpha:A_1\to A_0$, and the morphism $\eta:A\to
A^\vee[1]$ by an actual homomorphism of complexes $(\eta_1,\eta_0)$. Then
pick a point $P\in X$ and lift a basis of $\cok(\alpha)(P)$ to
$A_0$. replace $A_0$ by the free $\O_X$-module on this bases, and pull back
to get a quasi-isomorphic complex.  

Now any representative of $\eta$ has, necessarily, that $\eta_0$ is an
isomorphism in a neighborhood of $P$.  By Lemma~\ref{jfdsk}, we can assume
that $\eta_1=\eta_0^\vee$. Then both $\eta_0$ and $\eta_1$ are
isomorphisms of vector bundles.  Now use $\eta_0$ to identify $A_0$ with
$A_1^\vee$. 
\end{pf}

\subsection{Isometries}

\begin{defn}
Consider perfect complexes $A$ and $B$, endowed with non-degenerate
symmetric forms $\theta:A\to A^\vee[1]$ and $\eta:B\to B^\vee[1]$. An
isomorphism $\Phi:B\to A$, such that the diagram
$$\xymatrix{
 B\dto_{\eta} \rrto^\Phi && A\dto^\theta\\
 B^\vee[1] && A^\vee[1]\llto_{\Phi^\vee[1]}}$$
commutes in $D(X)$, is called an {\bf   isometry
 }$\Phi:(B,\eta)\to(A,\theta)$. 
\end{defn}

Note that because $\eta$ and $\theta$ are
isomorphisms, the condition on $\Phi$ is equivalent to
$\Phi^{-1}=\Phi^\vee[1]$, if we use $\eta$ and $\theta$ to identify
$A$ with $B$.

We include the following lemma on the local structure of
isometries for the information of the reader.  Since we do not use it
in the sequel, we omit the (lengthy) proof.

\begin{lem}
Let $A$ and $B$ be perfect, of amplitude contained in $[-1,0]$. 
Suppose $\theta:A\to A^\vee[1]$ and $\eta:B\to B^\vee[1]$ are
non-degenerate symmetric forms.  Let
$\Phi:B\to A$ be an isometry.

Suppose that $(A,\theta)$ and $(B,\eta)$ are represented as in
Example~\ref{simple} or Lemma~\ref{even.better}. 
Thus, $A=[E\stackrel{\alpha}{\to}E^\vee]$ and
$B=[F\stackrel{\beta}{\to}F^\vee]$, for vector bundles $E$ and $F$ on
$X$. Moreover, $\theta$ and $\eta$ are the respective identities.

Assume that that $\rk(F)=\rk(E)$.  
Then, \'etale locally in $X$ (Zariski locally if $X$ is a scheme), we
can find a vector bundle isomorphism
$$\phi:F\longto E\,,$$
such that $\alpha\circ\phi=\phi^{\vee-1}\circ\beta$, and
$(\phi,\phi^{\vee-1})$ represents $\Phi$:
$$\xymatrix{
B\dto_{\Phi}\ar@{}[r]|{=} &[F\rto^\beta\dto_\phi &
  F^\vee]\dto^{\phi^{\vee-1}}\\ 
A\ar@{}[r]|{=} & [E\rto^\alpha & E^\vee]}$$
In particular, $(\phi^{-1},\phi^\vee)$ represents $\Phi^\vee[1]$.
\end{lem}

\subsection{Symmetric obstruction theories}

Recall \cite{BF} that a {\em perfect obstruction theory }for the
scheme (or Deligne-Mumford stack) $X$ is a morphism $\phi:E\to L_X$ in
$D(X)$, where $E$ is perfect, of amplitude in $[-1,0]$, we have
$h^0(\phi):h^0(E)\to \Omega_X$ is an isomorphism and $h^{-1}:h^{-1}(E)\to
h^{-1}(L_X)$ is onto.

We denote the coherent sheaf $h^1(E^\vee)$ by $ob$ and call it the
{\em obstruction sheaf }of the obstruction theory. It contains in a
natural way the obstructions to the smoothness of $X$. Even though we
do not include $E$ in the notation, $ob$ is by no means an intrinsic
invariant of $X$. 

Any perfect obstruction theory for $X$ induces a virtual fundamental
class $[X]^\vir$ for $X$.  We leave the obstruction theory out of the
notation, even though $[X]^\vir$ depends on it. The virtual
fundamental class is an element of $A_{\rk E}(X)$, the Chow group of
algebraic cycles modulo rational equivalence. The degree of $[X]^\vir$
is equal to the rank of $E$.

\begin{defn} 
Let $X$ be a Deligne-Mumford stack. A {\em symmetric obstruction
  theory } for $X$ is a
triple $(E,\phi,\theta)$ where 
$\phi:E\to L_X$ is a perfect obstruction theory for $X$ and
$\theta:E\to  E^\vee[1]$ a non-degenerate symmetric bilinear form.
\end{defn}

We will often refer to such an $E$ as a symmetric obstruction theory,
leaving the morphisms $\phi$ and $\theta$ out of the notation.

\begin{rmk}
It is shown in \cite{Beh}, that for symmetric obstruction theories,
the virtual fundamental class  is intrinsic to $X$,
namely it is the degree zero Aluffi class of $X$.  
\end{rmk}

\begin{prop}
Every symmetric obstruction theory has expected dimension
zero. 
\end{prop}
\begin{pf}
Recall that the expected dimension of $E\to L_X$ is the rank of
$E$. If $E\to L_X$ is symmetric, we have $\rk E=\rk (E^\vee[1])=-\rk
E^\vee=-\rk E$ and hence $\rk E=0$.
\end{pf}

By this proposition, the following definition makes sense.

\begin{defn}
Assume $X$ is proper and  we have given a 
symmetric obstruction theory for $X$.
We define the {\em virtual count }of $X$ to be the number
$$\#^\vir(X)=\deg [X]^\mathrm{vir}=\int_{[X]^\mathrm{vir}}1\,.$$
If $X$ is a scheme (or an algebraic space), the virtual count is an
integer. In general it may be  a rational number.
\end{defn}

\begin{prop}
For a symmetric obstruction theory $E\to L_X$, the obstruction sheaf
is canonically isomorphic to the sheaf of differentials:
$$ob=\Omega_X\,.$$
\end{prop}
\begin{pf}
We have $ob=h^1(E^\vee)=h^0(E^\vee[1])=h^0(E)=\Omega_X$.
\end{pf}

\begin{cor}
For a symmetric obstruction theory we have
$h^{-1}(E)=\sheafhom(\Omega_X,\O_X)=T_X$.
\end{cor}
\begin{pf}
We always have that $h^{-1}(E)=ob^\vee$.
\end{pf}

\begin{defn}
Let $E$ and $F$ be symmetric obstruction theories for $X$. An {\em
  isomorphism }of symmetric obstruction theories is an isometry
  $\Phi:E\to F$ commuting with the maps to $L_X$.
\end{defn}

\begin{rem} 
Let $f:X\to X'$ be an \'etale morphism of Deligne-Mumford stacks, and
suppose that 
$X'$ has a symmetric obstruction theory $E'$. Then $f^\ast E'$ is
naturally a symmetric obstruction theory for $X$. 

Conversely, if we 
are given symmetric obstruction theories $E$ for $X$ and $E'$ for
$X'$, we will say that 
the morphism $f$ is compatible with the obstruction theories if $E$ is
isomorphic to $f^\ast E'$ as symmetric obstruction theory.
\end{rem}

\begin{rem} 
If $X$ and $X'$ are Deligne-Mumford stacks with symmetric obstruction
theories $E$ 
and $E'$, then $p_1^*E\oplus p_2^*E'$ is naturally a symmetric
obstruction theory for $X\times X'$. We call it the product
symmetric obstruction theory.
\end{rem}

\begin{ex}\Label{almclo}
Let $M$ be smooth and $\omega$ a closed 1-form on $M$.  Let
$X=Z(\omega)$ be the scheme-theoretic zero locus of $\omega$. 
Consider $\omega$ as a linear epimorphism $\omega^\vee:T_M\to I$,
where $I$ is the ideal sheaf of $X$ in $M$. Let us denote the
restriction to $X$ of the composition of $\omega^\vee$ and $d:I\to
\Omega_M$ by $\nabla\omega$. It is a linear homomorphism of vector
bundles $\nabla\omega:T_M|_X\to \Omega_M|_X$. Because $\omega$ is
closed, $\nabla\omega$ is symmetric and, as we have seen in
Example~\ref{simple}, defines a
symmetric bilinear form on the complex $E=[T_M|_X\to \Omega_M|_X]$. 

The morphism $\phi:E\to L_X$ as in the diagram
$$\xymatrix{
E\dto_\phi\ar@{}[r]|{=}& [T_M|_X \rto^-{\nabla\omega}\dto_{\omega^\vee}
& \Omega_M|_X]\dto^1\\ 
L_X\ar@{}[r]|{=} & [I/I^2\rto^-{d} & \Omega_M|_X]}$$
makes $E$ into a symmetric obstruction theory for $X$. In particular,
note that Example~\ref{hess} gives rise to a symmetric obstruction
theory on the Jacobian locus of a regular function. 

Let us remark that for the symmetry of $\nabla\omega$ and hence the
symmetry of the obstruction theory given by $\omega$, it is sufficient
that $\omega$ be {\em almost closed}, which means that $d\omega \in
I\Omega^2_M$. 
\end{ex}

\subsection{A remark on the lci case}

We will show that the existence of a symmetric obstruction theory puts
strong restrictions on the singularities $X$ can have.

For the following proposition, it is important to recall that we are
working in characteristic zero. 

\begin{prop}
Let $E\to L_X$ be a perfect obstruction theory, symmetric or not.  A
criterion for the obstruction sheaf to be locally free is that
$X$ be a reduced local complete intersection. 
\end{prop}
\begin{pf}
As the claim is local, we may assume that $E$ has a global resolution
$E=[E_1\to E_0]$, that  $X\hookrightarrow M$ is embedded in a smooth
scheme $M$ (with ideal sheaf $I$) and that $E\to L_X$ is given by a
homomorphism of complexes $[E_1\to E_0]\longrightarrow
[I/I^2\to\Omega_M|_X]$. We may even assume that $E_0\to\Omega_M|_X$ is
an isomorphism of vector bundles. 

Under the assumption that $X$ is a reduced local complete
intersection, we have that $I/I^2$ is locally free and that $I/I^2\to
\Omega_M|_X$ is injective. Then a simple diagram chase proves that we
have a short exact sequence of coherent sheaves
$$0\longrightarrow h^{-1}(E)\longrightarrow E_1\longrightarrow
I/I^2\longrightarrow 0\,.$$
Hence, $h^{-1}(E)$ is a subbundle of $E_1$ and $ob=h^{-1}(E)^\vee$. In
particular, $ob$ is locally free. 

{\em Remark. }We always have that $h^{-1}(E)=ob^\vee$, the converse
is generally false.
\end{pf}

\begin{cor}
If $X$ is  a reduced local complete intersection and admits a
symmetric obstruction theory, then $X$ is smooth.
\end{cor}
\begin{pf}
Because $ob=\Omega_X$, we have that $\Omega_X$ is locally free.  This
implies that $X$ is smooth.
\end{pf}

\subsection{Examples}

\subsubsection{Lagrangian intersections}

Let $M$ be an algebraic  symplectic manifold and $V$, $W$ two Lagrangian
submanifolds.  Let $X$ be the scheme-theoretic intersection. Then $X$
carries a canonical symmetric obstruction theory. 

To see this, note first of all that for a Lagrangian submanifold
$V\subset M$,  the normal bundle is equal to the
cotangent bundle, $N_{V/M}=\Omega_V$. The isomorphism is given
by $v \longmapsto \sigma(v,-)$, where $\sigma$ is the symplectic form,
which maps $N_{V/M}=T_M/T_V$ to $\Omega_V=T_V^\vee$. It is essentially
the definition of Lagrangian, that this map is an isomorphism of
vector bundles on $V$.

Next, note that the obstruction theory for $X$ as an intersection of
$V$ and $W$ can be represented by the complex
$$\xymatrix@1{E\ar@{}[r]|-=&[\,\Omega_M\rrto^-{\res_V-\res_W}&&
    \Omega_V\oplus\Omega_W\,]\big|_X\,.}$$  
The shifted dual is
$$\xymatrix@1{E^\vee[1]\ar@{}[r]|-=&[\,T_V\oplus T_W\rto&
    T_M\,]\big|_X\,.}$$ 
Define $\theta:T_M\to \Omega_V\oplus \Omega_W$ as the  canonical map
    $T_M\to N_{V/M}\oplus N_{W/M}$ given by the 
projections, multiplied by the scalar factor $\frac{1}{2}$. 
Then $(\theta^\vee,\theta)$ defines a morphism of
complexes $E^\vee[1]\to E^\vee$:
$$\xymatrix{
E^\vee[1]\ar@{}[r]|=\dto & {[T_V\oplus T_W}\rto\dto_{\theta^\vee} &
  {\phantom{M}T_M\phantom{M}]} \dto^\theta\\
E\ar@{}[r]|= & {[\phantom{M}\Omega_M\phantom{M}}\rto &
  {\Omega_V\oplus\Omega_W]}}$$ 
One checks that $(\theta^\vee,\theta)$ is a quasi-isomorphism. 

As $(\theta^\vee,\theta)^\vee[1]=(\theta^\vee,\theta)$, this morphism
of complexes defines a symmetric bilinear form on $E^\vee[1]$, hence
on $E$. Thus  $E$ is a symmetric obstruction theory on $X$.

\subsubsection{Sheaves on Calabi-Yau threefolds}

Let $Y$ be an integral proper 3-dimensional  Gorenstein
Deligne-Mumford stack (for example a projective threefold). 
By the Gorenstein assumption, $Y$ admits a dualizing sheaf $\omega_Y$,
which is a line bundle over $Y$, also called the {\em canonical
}bundle.  Let $\omega_Y\to \O_Y$ be a non-zero homomorphism, giving
rise the the short exact sequence
$$\xymatrix@1{0\rto &  \omega_Y\rto & \O_Y\rto & \O_D\rto & 0}\,,$$
so that $D$ is an anti-canonical divisor on $Y$. In fact, $D$ is a
Cartier divisor.  Of course, $D$ may be empty (this case we refer to
as the Calabi-Yau case). 

Finally, choose an arbitrary line bundle $L$ on $Y$.  Often we are
only interested in the case $L=\O_Y$. 

Now let us define a certain moduli stack $\MM$  of sheaves on $Y$.
For an 
arbitrary $\cc$-scheme $S$, let  $\MM(S)$ be the groupoid of pairs
$(\mathcal{E},\phi)$. Here 
$\mathcal{E}$ is a sheaf of $\O_{Y\times S}$-modules, such that

(i) $\mathcal{E}$ coherent,

(ii) $\mathcal{E}$ is flat over $S$,

(iii) $\mathcal{E}$ is perfect as an object of the derived category
of $Y\times S$, i.e., locally admits finite free resolutions, (by
Cor.~4.6.1 of Exp.~III of SGA~6, this is a condition which may be
checked on the fibres of $\pi:Y\times S\to S$). 

The second component of the pair $(\mathcal{E},\phi)$ is an
isomorphism  $\phi:\det\mathcal{E}\to L$  of line
bundles on 
$Y\times S$.  Note that the determinant $\det\mathcal{E}$ is
well-defined, by Condition~(iii) on $\mathcal{E}$. 

We require two more conditions on $\mathcal{E}$, namely that for every
point $s\in S$, denoting the fibre of $\mathcal{E}$ over $s$ by
$\mathcal{E}_s$, we have

(iv) $\mathcal{E}_s$ is simple, i.e., $\kappa(s)\to
\Hom(\mathcal{E}_s,\mathcal{E}_s)$ is an isomorphism,

(v) the map induced by the trace
$R\sheafhom(\mathcal{E}_s,\mathcal{E}_s)\to\O_{Y_s}$ is an isomorphism in a
neighborhood of $D_s$. 

The last condition (v) is empty in the Calabi-Yau case. It is, for
example, satisfied if $\mathcal{E}_s$ is locally free of rank 1 in a
neighborhood of $D$. 

We let $X$ be an open substack of $\MM$ which is algebraic 
(for example, fix the Hilbert polynomial and  pass to stable objects, but
we do not want to get more restrictive than necessary).  Then $X$ is a
Deligne-Mumford stack.   We will now construct a symmetric obstruction
theory for $X$.  

For this, denote the universal sheaf on $Y\times X$ by
$\mathcal{E}$ and the projection $Y\times X\to X$ by $\pi$. 
Consider the trace map 
$R\sheafhom(\mathcal{E},\mathcal{E})\to \O$ and let $\mathcal{F}$ be
its shifted cone, so that we obtain a distinguished triangle in
$D(\O_{Y\times X})$:
$$\xymatrix@C=1pc{
& \O\dlto_{+1} & \\
\mathcal{F}\rrto && R\sheafhom(\mathcal{E},\mathcal{E})\ulto_{tr}}$$
Note that $\mathcal{F}$ is self-dual: $\mathcal{F}^\vee=\mathcal{F}$,
canonically. 

\begin{lem}
The complex
$$E=R\pi_\ast R\sheafhom(\mathcal{F},\omega_Y)[2]$$
is an obstruction theory for $X$.
\end{lem}
\begin{pf}
This is well-known deformation theory. See, for example, \cite{RT}.
\end{pf}

The homomorphism $\omega_Y\to\O_Y$ induces an {\em isomorphism }
$$R\pi_\ast(\mathcal{F}\otimes \omega_Y)\longrightarrow
R\pi_\ast\mathcal{F}\,,$$
because the cone if this homomorphism is
$R\pi_\ast(\mathcal{F}\otimes\O_D)$ and $\mathcal{F}\otimes\O_D=0$, by
Assumption~(v), above. Dualizing and shifting, we obtain an isomorphism
$$(R\pi_\ast\mathcal{F})^\vee[-1]\longrightarrow
\big(R\pi_\ast(\mathcal{F}\otimes \omega_Y)\big)^\vee[-1]\,.$$
Exploiting the fact that $\mathcal{F}$ is self-dual, we may rewrite
this as
$$(R\pi_\ast\mathcal{F})^\vee[-1]\longrightarrow
\big(R\pi_\ast R\sheafhom(\mathcal{F},\omega_Y)\big)^\vee[-1]\,,$$
or in  other words
\begin{equation}\Label{iow}
(R\pi_\ast\mathcal{F})^\vee[-1]\longrightarrow E^\vee[1]\,.
\end{equation}

Now, relative Serre duality
for the morphism $\pi:Y\times X\to X$
applied to $\mathcal{F}$ states that
$$R\pi_\ast
R\sheafhom(\mathcal{F},\omega_Y[3])=(R\pi_\ast\mathcal{F})^\vee\,,$$ 
or in other words
$$E=(R\pi_\ast\mathcal{F})^\vee[-1]\,.$$
Thus, we see that (\ref{iow}) is, in fact, an isomorphism
$$\theta:E\longrightarrow E^\vee[1]\,.$$

\begin{lem}
The isomorphism $\theta:E\to E^\vee[1]$ satisfies the symmetry
property $\theta^\vee[1]=\theta$. 
\end{lem}
\begin{pf}
This is just a derived version of the well-known fact that
$tr(AB)=tr(BA)$, for endomorphisms $A$, $B$ of a free module.
\end{pf}

\begin{lem}
The complex $E$ has perfect amplitude contained in the  interval
$[-1,0]$.
\end{lem}
\begin{pf}
Perfection is clear.  To check the interval, note that by symmetry of
$E$ it suffices to check that the interval is $[-1,\infty]$. We have
seen that $E=R\pi_\ast\mathcal{F}[2]$.  So the interval is no wore
than $[-2,\infty]$.  But $h^{-2}(E)=0$, by Assumption~(iv), above.
\end{pf}

\begin{cor}\Label{sot}
The Deligne-Mumford stack $X$  admits, in a natural way, a
symmetric obstruction theory, namely 
$$E=R\pi_\ast
R\sheafhom(\mathcal{F},\omega_Y)[2]=R\pi_\ast\mathcal{F}[2]\,.
$$
\end{cor}

We call this obstruction theory the {\em Donaldson-Thomas }obstruction
theory. 

In the next two propositions we mention two special cases.  The first
was originally introduced in \cite{RT}, where the symmetry was pointed
out, too.

\begin{prop}
Let $Y$ be a smooth projective threefold with trivial canonical
bundle, and let $X$ be the fine moduli stack of stable sheaves on $Y$
of rank $r>0$, with fixed determinant $L$ and with Chern classes
$c_2$, $c_3$. Then $X$ admits a symmetric obstruction theory.
\end{prop}
\begin{pf}
In fact, every trivialization $\omega_Y=\O_Y$ induces a symmetric
obstruction theory. 
\end{pf}

\begin{prop}\Label{dtot}
Let $Y$ be a smooth projective threefold and $D$ an effective
anti\-canonical divisor 
on $Y$. Let $X'$ be the scheme of torsion-free rank 1 sheaves with
trivial determinant and fixed Chern classes $c_2$, $c_3$. 
Recall that such sheaves can be identified with ideal sheaves. Let
$X\subset X'$ be the open subscheme consisting of ideal sheaves which
define a subscheme of $Y$ whose support is disjoint from $D$. 
Then $X$ admits a symmetric obstruction theory. 

For example, $\Hilb^n(Y\setminus D)$, the Hilbert scheme of length $n$
subschemes of $Y\setminus D$ admits a symmetric obstruction theory.
\end{prop}
\begin{pf}
Again, we would like to point out that every homomorphism $\omega_Y\to
\O_Y$ defining $D$ gives rise to a symmetric obstruction theory on
$X$.  Even though the compactification is used in
its construction, this symmetric obstruction theory does not depend on
which compactification is chosen. 
\end{pf}

\subsubsection{Stable maps to Calabi-Yau threefolds}

\begin{prop} 
Let $Y$ be a Calabi-Yau threefold and let $X$ be the open locus in the 
moduli stack of stable maps parameterizing immersions of smooth
curves. Then the Gromov-Witten obstruction theory of $X$ is symmetric,
in a natural way.
\end{prop}
\begin{pf}
Let $\pi:C\to X$ be the universal curve and $f:C\to Y$ the universal
map. Let $F$ be the kernel of $f^\ast\Omega_Y\to \Omega_C$, which is a
vector bundle of rank 2 on $C$. The Gromov-Witten obstruction theory
on $X$ is $E=R\pi_\ast(F\otimes\omega_{C/X})[1]$. By Serre duality for
$\pi:C\to X$, we have $E^\vee[1]=R\pi_\ast(F^\vee)[1]$. 

As $F$ is of rank 2, we have $F=F^\vee\otimes \det F$. Because $Y$ is
Calabi-Yau, we have $\det F\otimes \omega_{C/X}=\O_C$. Putting these
two facts together, we get $F\otimes\omega_{C/X}=F^\vee$ and hence
$E=E^\vee[1]$. 
\end{pf}

\eject
\section{Equivariant symmetric obstruction theories}

\subsection{A few remarks on equivariant derived categories}

Let $X$ be a scheme with an action of an algebraic group $G$.  Let
$(\qcoh\O_X)^G$ denote the abelian category of $G$-equivariant
quasi-coherent $\O_X$-modules. Thus, and object of $(\qcoh\O_X)^G$ is
a quasi-coherent $\O_X$-module $F$ together with descent data to the
quotient stack $[X/G]$, in other words and isomorphism between $p^\ast
F$ and $\sigma^\ast F$ satisfying the cocycle condition. Here $p$ and
$\sigma$ are projection and action maps $X\times G\to X$,
respectively.  Denote by $D(X)^G$ the derived category of
$(\qcoh\O_X)^G$. Note that $\O_X$ is an object of $D(X)^G$, in a
natural way.

There is the forgetful functor $D(X)^G\to D(X)$, which maps a complex
of equivariant sheaves to its underlying complex of sheaves. It is an
exact functor.

To simplify matters, let us make two assumptions:

\noindent (a) $X$ admits a $G$-equivariant ample invertible sheaf
$\O(1)$,

\noindent (b) $G$ is a diagonalizable group,
i.e.,  $G=\spec \cc[W]$ is the spectrum of the group ring of a
finitely generated abelian group $W$.  Then $W$ is canonically
identified with the character group of $G$.

\subsubsection{The affine case}

If $X=\spec A$ is  affine,  $A$ is $W$-graded. A
$G$-equivariant $\O_X$-module is the same thing as a $W$-graded
$A$-module. 

We call a $W$-graded $A$-module {\em quasi-free}, if it is
free as an $A$-module on a set of homogeneous generators. Any
quasi-free $W$-graded $A$-module is isomorphic to a direct sum of
shifted copies of $A$. Quasi-free $W$-graded $A$-modules are
projective objects in the abelian category $(\qcoh\O_X)^G$ of
$W$-graded $A$-modules.  Hence this category has enough projective
objects.

\subsubsection{The global case}

Let $F$ be a $G$-equivariant $\O_X$-module. We can shift $F$ by any
character $w\in W$ of $G$. We denote the shift by $F[w]$. Every
$G$-equivariant quasi-coherent $\O_X$-module $F$ can be written as a
quotient of sheaf of the form 
\begin{equation}\Label{form}
\bigoplus_{i\in I}\O(n_i)[w_i]\,.
\end{equation}
Thus, every $G$-equivariant quasi-coherent $\O_X$-module admits
left resolutions consisting of objects of form (\ref{form}). More
generally, every bounded above complex in $D(X)^G$ can be replaced by
a bounded above complex of objects of type (\ref{form}). These
resolutions are $G$-equivariant. 

Since objects of the form (\ref{form}) are locally free as
$\O_X$-modules (forgetting the $G$-structure), we can use these
resolutions to compute the derived functors of $\otimes$ and 
$\sheafhom(-,F)$. Thus we see that for $G$-equivariant quasi-coherent
$\O_X$-modules  $E$, $F$ the quasi-coherent $\O_X$-modules
$\sheaftor_i(E,F)$ and $\sheafext^i(E,F)$ are again
$G$-equivariant. More generally, we see that for bounded above objects
$E$, $F$ of $D(X)^G$, the objects $E\ltensor F$ and $R\sheafhom(E,F)$
are again in $D(X)^G$. 

For a $G$-equivariant sheaf $E$, we write $E^\vee=\sheafhom(E,\O_X)$.
For a bounded above object $E$ of $D(X)^G$, we write
$E^\vee=R\sheafhom(E,\O_X)$.

Let $\{U_i\}$ be an invariant affine open cover. Let $F$ be a
$G$-equivariant quasi-coherent $\O_X$-module. Then, the \v Cech
resolution  $\cC^\bullet(\{U_i\},F)$ is a right resolution of $F$ by
$G$-equivariant quasi-coherent $\O_X$-modules. It is an acyclic
resolution for the global section functor, showing that the cohomology
groups $H^i(X,F)$ are $W$-graded. 
More generally, let $f:X\to Y$ be a $G$-equivariant morphism. Then we
see that $R^if_\ast F$ are $G$-equivariant quasi-coherent
$\O_Y$-modules. 

Moreover, if $E$ is a bounded below complex in $D(X)^G$, we can
construct the associated \v Cech complex $\cC(\{U_i\},E)$, which is a
double complex.  Passing to the associated single complex, we see that
we may replace $E$ by a bounded below complex of $G$-equivariant
$\O_X$-modules which are acyclic for $f_\ast$, for any $G$-equivariant
morphism $f:X\to Y$. Thus we see that the functor $Rf:D(X)\to D(Y)$
passes to a functor $Rf:D(X)^G\to D(Y)^G$. 

\subsubsection{The cotangent complex}

If $X$ is a $G$-scheme as above, the sheaf of K\"ahler differentials
$\Omega_X$ and its dual $T_X=\Omega_X^\vee$ are $G$-equivariant. 

We can use the equivariant ample line bundle $L$ to construct a
$G$-equivariant embedding $X\hookrightarrow M$ into a smooth
$G$-scheme $M$. The cotangent complex $I/I^2\to\Omega_M|X$ is then
naturally an object of $D(X)^G$. The usual proof that $L_X$ is a
canonically defined object of $D(X)$ works equivariantly and proves
that $L_X$ is a canonically defined object of $D(X)^G$.  By
canonically defined, we mean that any two constructions are related by
a canonical isomorphism.

\subsubsection{Perfect objects}

We call an object $E$ in $D(X)^G$ {\em perfect\/} (of perfect
amplitude in the interval $[m,n]$), if its underlying
object of $D(X)$ is perfect (of perfect amplitude in the interval
$[m,n]$). 

\begin{rmk}
If $X$ is a scheme and $E$ in $D(X)$ is a perfect complex, of perfect
amplitude contained in $[m,n]$, then we can write $E$ locally as a
complex $[E^m\to\ldots\to E^n]$ of free $\O_X$-modules contained in
the interval $[m,n]$.  This is essentially because if $E\to E''$ is an
epimorphism of locally free coherent sheaves, the kernel is again
locally free coherent.

In the equivariant context, we have to forgo this convenient fact.
Suppose $E$ in $D(X)^G$ is perfect, again of amplitude contained in
$[m,n]$. We can, as we saw above, write $E$ as a bounded above complex
of sheaves of form (\ref{form}), all of them coherent, i.e., with
finite indexing set $I$.  But when we cut off
this infinite complex to fit into the interval $[m,n]$, we end up with
a $G$-equivariant quasi-coherent sheaf which is locally free coherent
as an $\O_X$-module without the $G$-structure, but which is not
locally quasi-free and not locally projective in the category
$(\qcoh\O_X)^G$. 
\end{rmk}

\subsection{Symmetric equivariant obstruction theories}

\begin{defn}
Let $X$ be a scheme with  a $G$-action. An {\em
equivariant }perfect obstruction theory is a morphism $E\to L_X$ in
the category $D(X)^G$, which is a perfect obstruction theory as a
morphism in $D(X)$. (This definition is originally due to
Graber-Pandharipande~\cite{grp}.)

A {\em symmetric equivariant }obstruction theory, (or an {\em
equivariant symmetric }obstruction theory) is a pair $(E\to L_X,E\to
E^\vee[1])$ of morphisms in the category $D(X)^G$, such that $E\to
L_X$ is an (equivariant) perfect obstruction theory and $\theta:E\to
E^\vee$ is an isomorphism satisfying $\theta^\vee[1]=\theta$. 
\end{defn}

Note that this is more than requiring that the obstruction theory be
equivariant and symmetric, separately, as we can see in the following
example.

\begin{ex}\Label{almost}
Let $\omega=\sum_i^n f_i dx_i$ be an almost closed 1-form on
$\A^n$. Recall from Example~\ref{almclo} that $\omega$ defines a
symmetric obstruction theory 
$H(\omega)=[T_M|_X\stackrel{\nabla\omega}{\longrightarrow}
\Omega_M|_X]$ on the zero locus $X$ of $\omega$.

Define a $\Gm$-action on $\A^n$ by setting the degree of $x_i$
to be $r_i$, where $r_i\in\zz$. Assume that each $f_i$ is homogeneous
with respect to these degrees and denote the degree of $f_i$ by $n_i$. 
Then the zero locus $X$ of $\omega$ inherits a $\Gm$-action.

If we let $\Gm$ act on $T_M$ by
declaring the degree of $\frac{\del}{\del x_i}$ to be equal to $n_i$,
then $H(\omega)$ is $\Gm$-equivariant as well as the morphism
$H(\omega)\to L_X$.  Thus $H(\omega)$ is an equivariant obstruction
theory. 

But note that $H(\omega)$ is not equivariant symmetric. This is because
the identity on $H(\omega)$ (which is $\theta$ in this case) is not
$\Gm$-equivariant if we consider it as a homomorphism $H(\omega)\to
H(\omega)^\vee[1]$. Unless $n_i=-r_i$, because then the degree of
$\frac{\del}{\del x_i}$ is equal to its degree as the dual of $dx_i$. 

In the case $n_i=-r_i$, the form $\omega=\sum f_i dx_i$ is an
invariant element of $\Gamma(M,\Omega_M)$, or an equivariant
homomorphism $\O_M\to \Omega_M$. In this case we do get an equivariant
symmetric obstruction theory.
\end{ex}

\subsubsection{The equivariant Donaldson-Thomas obstruction theory}

Let $G$ be a diagonalizable group as above.  Consider a projective
threefold $Y$, endowed with a linear $G$-action.  Consider a
$G$-equivariant non-zero homomorphism $\omega_Y\to \O_Y$, defining the
$G$-invariant anti-canonical Cartier divisor $D$.

\begin{prop}\Label{edt}
Let $X$ be as in Proposition~\ref{dtot}. Then the Donaldson-Thomas
obstruction theory of Corollary~\ref{sot} on $X$ is $G$-equivariant
symmetric.
\end{prop}
\begin{proof} 
Let $X'$ be the compactification of $X$ as in Proposition~\ref{dtot}.
Let $\eE$ be the universal sheaf on $Y\times X$ and $Z\subset Y\times
X$ be the universal subscheme.  We have an exact sequence 
$$\xymatrix@1{0\rto & \eE\rto & \O_{Y\times X}\rto & \O_Z\rto &
  0}\,.$$
Let $\pi:Y\times X\to X$ be the projection.  Note
that $\eE$ and $\O_Z$ are $G$-equivariant. This follows directly from the
universal mapping property of $\eE$.  

The standard ample invertible sheaf on $X'$ is $\det\pi_\ast(\O_Z(n))$,
for $n$ sufficiently large.  It is $G$-equivariant, as all ingredients
in its construction are.  Hence $X$ admits an equivariant ample
invertible sheaf.  

Next, notice that all the constructions involved in producing the
obstruction theory $E=R\pi_\ast R\sheafhom(\fF,\omega_Y)[2]$ work
equivariantly. Hence the symmetric obstruction theory is equivariant.

To prove that it is equivariant symmetric, we just need to remark that
the bilinear form $\theta$ is induced from $\omega\to \O_Y$, which is
equivariant, and that Serre duality is equivariant, because it is
natural. 
\end{proof}

\subsection{Local structure in the $\Gm$-case}

Let $G=\Gm$.  We will prove that Example~\ref{almost} describes the
unique example of a symmetric $\Gm$-equivariant obstruction theory, at
least locally around a fixed point.

\begin{lem}\Label{preparation}
Let $X$ be an affine $\Gm$-scheme with a fixed point $P$. Let $n$
denote the dimension of $T_X|P$, the Zariski tangent space of $X$ at
$P$. Then there
exists an invariant affine open neighborhood $X'$ of $P$ in $X$, a
smooth $\Gm$-scheme $M$ of dimension $n$ and an equivariant closed
embedding $X'\hookrightarrow M$
\end{lem}
\begin{pf}
Let $A$ be the affine coordinate ring of $X$.  The $\Gm$-action
induces a grading on $A$. Let $\Mm$ be the maximal ideal given by the
point $P$. We can lift an eigenbasis of $\Mm/\Mm^2$ to homogeneous
elements $x_1,\ldots,x_n$ of $\Mm$. Choose homogeneous elements
$y_1,\ldots,y_m$  in $\Mm$ in such a way that
$x_1,\ldots,x_n,y_1,\ldots,y_m$ is a set of generators of $A$ as a
$\cc$-algebra. This defines a closed embedding $X\hookrightarrow
\A^{n+m}$, which is equivariant if we define a $\Gm$-action on
$A^{n+m}$ in a suitable, obvious, way. 

We have thus written $A$ as a quotient of $\cc[x,y]$.  Let $I$ denote
the corresponding homogeneous ideal in $\cc[x,y]$. Then we have
$$\Mm/\Mm^2=(x,y)/\big(I+(x,y)^2\big)\,.$$
Since this $\cc$-vector space is generated by $x_1,\ldots,x_n$, we
have, in fact:
$$y_i\in I+(x,y)^2+(x)\,,$$
for $i=1,\ldots,m$. We can therefore find homogeneous elements
$f_1,\ldots,f_m\in I$, such that 
$$y_i-f_i\in (x,y)^2+(x)\quad\text{and}\quad\deg f_i=\deg y_i\,,$$
for all $i=1,\ldots,m$. Let $g\in \cc[x,y]$ be the determinant of the
Jacobian matrix $(\frac{\del f_i}{\del y_j})$.  We see that $g$ is
homogeneous of degree $0$ and that 
$g(0,0)=1$. Let $U\subset \A^{n+m}$ be the affine open subset where
$g$ does not vanish.  This is an invariant subset containing $P$. Let
$Z\subset \A^{n+m}$ be the closed subscheme defined by
$(f_1,\ldots,f_m)$.  It carries an induced $\Gm$-action. The
intersection $M=Z\cap U$ is a smooth scheme of dimension $n$. 

As $(f_1,\ldots,f_m)\subset I$, we have that $X$ is a closed subscheme
of $Z$. Let $X'=X\cap U$. 
\end{pf}

\begin{prop}\Label{annoying}
Let $X$ be an affine $\Gm$-scheme with a fixed point $P$ and let
$n=\dim T_X|_P$. Furthermore, let $X$ be endowed with a symmetric
equivariant obstruction theory $E\to L_X$. Then there exists an
invariant affine open neighborhood $X'$ of $P$ in $X$, an equivariant
closed embedding $X'\hookrightarrow M$ into a smooth $\Gm$-scheme $M$
of dimension $n$ and  an invariant almost closed 1-form $\omega$ on $M$
such that $X=Z(\omega)$.  We can further construct an equivariant
isometry $E\to H(\omega)$ commuting with the maps to $L_X$, but it
will not be necessary for the purposes of this paper. 
\end{prop}
\begin{pf}
We apply Lemma~\ref{preparation}, to obtain the equivariant closed
embedding $X'\hookrightarrow M$.  Write $X$ for $X'$. 
Let $A$ be the affine coordinate ring of $X$ and $I$ the
ideal of $\Gamma(\O_M)$ defining $X$.  

Consider the object $E$ of $D(X)^\Gm$. 
We  can represent $E$ by an infinite
complex $[\ldots\to E_1\to E_0]$ of finitely generated quasi-free
$A$-modules.

Because quasi-free modules are projective, if $E$ is represented by a
bounded above complex of quasi-free modules 
and $E\to F$ is a a morphism in $D(X)^G$, then $E\to F$ can be
represented by an actual morphism of complexes, without changing $E$. 

Thus we have morphisms of complexes of graded modules $[\ldots\to
E_1\to E_0]\to[I/I^2\to \Omega_M|_X]$ and
$\theta:[\ldots\to E_1\to  E_0]\to[E_0^\vee\to E_1^\vee\to\ldots]$.
We can represent the equality of derived category morphisms
$\theta^\vee[1]=\theta$ by a homotopy between $\theta^\vee[1]$ and
$\theta$, because $E$ is a bounded above complex of quasi-frees. Then,
as in the proof of Lemma~\ref{jfdsk}, we can replace $\theta_0$ by
$\frac{1}{2}(\theta_0+\theta_1^\vee)$ and $\theta_1$ by
$\frac{1}{2}(\theta_1+\theta_0^\vee)$, without changing the homotopy
class of $\theta$. Then we have that $\theta_1=\theta_0^\vee$. 

Now we can replace $E_1$ by $\cok(E_2\to E_1)$ and $E_1^\vee$ by
$\ker(E_1^\vee\to E_2^\vee)$.  Because of the perfection of $E$, both
$\cok(E_2\to E_1)$ and $\ker(E_1^\vee\to E_2^\vee)$ are projective
$A$-modules (after forgetting the grading), which are, moreover, dual
to each other.  

Thus we have now represented $E$ by a complex $[E_1\to E_0]$ of
equivariant vector bundles and $E\to L_X$ and $\theta:E\to E^\vee[1]$
by equivariant morphisms of complexes. Moreover,
$\theta=(\theta_0^\vee,\theta_0)$, for an equivariant morphism of
vector bundles $\theta_0:E_0\to E_1^\vee$. 

Now we remark that we may assume that the rank of $E_0$ is equal to
$n$. Simply lift a homogeneous basis of $\Omega_X|_P$ to $E_0$ and
replace $E_0$ by the quasi-free module on these $n$ elements of
$E_0$. Then pass to an invariant open neighborhood of $P$ over which
both $E_0\to \Omega_M|_X$ and $\theta_0:E_0\to E_1^\vee$ are
isomorphisms.  Use these isomorphisms to identify. Then our
obstruction theory is given by an equivariant homomorphism
$$\xymatrix{[T_M|_X\rto^\alpha\dto_\phi & \Omega_M|_X]\dto^{\id}\\
[I/I^2\rto &\Omega_M|_X]}$$
such that $\alpha^\vee=\alpha$. Note that $\phi$ is necessarily
surjective. 

As we may assume that $\Omega_M|_X$ and hence $T_M|_X$ is given by a
quasi-free $A$-module, we may lift $\phi$ to an equivariant
epimorphism $T_X\to I$.  This gives the invariant 1-form
$\omega$. 
\end{pf}

\eject
\section{The main theorem}

\subsection{Preliminaries on linking numbers}

Here our dimensions are all {\em real }dimensions.

We work with orbifolds.  Orbifolds are differentiable stacks of
Deligne-Mumford type, which means that they are representable by Lie
groupoids $X_1\toto X_0$, where source and target maps $X_1\to X_0$
are \'etale (i.e., local diffeomorphisms) and the diagonal $X_1\to
X_0\times X_0$ is proper. If a compact Lie group $G$ acts with finite
stabilizers on a manifold $X$, the quotient stack $[X/G]$ is an
orbifold.

All our orbifolds will tacitly assumed to be {\em oriented}, which
means that any presenting groupoid $X_1\toto X_0$ is oriented, i.e.,
$X_0$ and $X_1$ are oriented and all structure maps (in particular
source and target $X_1\to X_0$) preserve orientations. 

Given an orbifold $X$, presented by the groupoid $X_1\toto X_0$, with
proper diagonal $X_1\to X_0\times X_0$, the image of the diagonal is a
closed equivalence relation on $X_0$. The quotient is the {\em coarse
moduli space }of $X$.

We call an orbifold {\em compact}, if its course moduli space is
compact. More generally, we call a morphism $f:X\to Y$ of orbifolds
{\em proper}, if the induced map on coarse moduli spaces is proper. 

To fix ideas, let $H^\ast(X)$ denote de Rham cohomology of the
orbifold $X$.  For the definition and basic properties of this cohomology
theory, see \cite{trieste}.  Note that homotopy invariance holds:
the projection $X\times \rr\to X$ induces an
isomorphism $H^\ast(X)\to H^\ast(X\times \rr)$.

If $f:X\to Y$ is a proper morphism of orbifolds, there exists a
wrong way map $f_!:H^i(X)\to H^{i-d}(Y)$, where $d=\dim
X-\dim Y$ is the relative dimension of $f$.  If $Y$ is the point, then
we also denote $f_!$ by $\int_X$. We will need the following
properties of $f_!$:

(i) Functoriality: $(g\circ f)_!=g_!\circ f_!$. 

(ii) Naturality: if $v:V\subset Y$ is an open suborbifold and
$u:U\subset X$ the inverse image of $U$ under $f:X\to Y$, we have
$v^\ast\comp f_!=g_!\comp u^\ast$, where $g:U\to V$ is the restriction
of $f$.

(iii) Projection formula: $f_!\big(f^\ast(\alpha)\cup\beta\big)=\alpha\cup
f_!(\beta)$. 

(iv) Poincar\'e duality: if $X$ is a compact orbifold, the pairing
$\int_X\alpha\cup\beta$ between $H^i(X)$ and $H^{n-i}(X)$ is a perfect
pairing of finite dimensional $\rr$-vector spaces ($n=\dim X$). 

(v) Long exact sequence: if $\iota:Z\subset X$ is a closed suborbifold with
open complement $U$, there is a long exact sequence ($c=\dim X-\dim Z$)
$$\xymatrix@1{
\ldots\rto^-\del & \rto H^{i-c}(Z)\rto^{\iota_!} & H^i(X)\rto &
H^i(U)\rto^-\del & 
H^{i-c+1}(Z)\rto & \ldots}
$$
In the situation of (v), we call $\cl(Z)=\iota_!(1)\in H^c(X)$ the
{\em class }of $Z$. 

We could use any other cohomology theory with characteristic zero
coefficients which satisfies these basic properties.

\begin{rmk}\Label{diffhom}
Let $T\subset \rr$ be an open interval containing the points $0$ and
$1$. Let $Z$ and $X$ be a compact orbifolds and $h:Z\times T\to X$ a
differentiable 
morphism of orbifolds such that $h_0:Z\times\{0\}\to X$ and
$h_1:Z\times \{1\}\to X$ are isomorphisms onto closed suborbifolds
$Z_0$ and $Z_1$ of $X$. We call $h$ a {\em differentiable homotopy
}between $Z_0$ and $Z_1$.  It is not difficult to see, using
Poincar\'e duality and homotopy invariance, that the existence of such
an $h$ implies that $\cl(Z_0)=\cl(Z_1)\in H^\ast(X)$.
\end{rmk}

\subsubsection{Linking numbers and $S^1$-actions}

Let $A$ and $B$ be closed submanifolds, both of dimension $p$, of a
compact manifold $S$ of dimension $2p+1$.  Assume that
$H^{p+1}(S)=H^{p}(S)=0$ and that $A\cap B=\varnothing$. For
simplicity, assume also that $p$ is odd.  

Under these assumptions we can define the 
{\em linking number }$L_S(A,B)$ as follows. 
By our assumption, the boundary map $\del:H^{p}(S\setminus B)\to H^0(B)$ is
an isomorphism.  Let $\beta\in H^{p}(S\setminus B)$ the the  unique
element such that $\del\beta=1\in H^0(B)$. Via the inclusion $A\to
S\setminus B$ we restrict $\beta$ to $A$ and set
$$L_S(A,B)=\int_A \beta\,.$$

Now assume $A'$ is another closed submanifold of $S$ of dimension $p$,
and $A'\cap B=\varnothing$, too.  Thus $L_S(A',B)$ is defined.  We
wish to compare $L_S(A',B)$ with $L_S(A,B)$. 

Suppose $h:Z\times T\to S$ is a differentiable homotopy between $A$
and $A'$, as in Remark~\ref{diffhom}. It is an obvious, well-known
fact, that if the 
image of $h$ is entirely 
contained in $S\setminus B$, then $L_S(A',B)=L_S(A,B)$.  We wish to show that
in the presence of an $S^1$-action, 
$L_S(A',B)=L_S(A,B)$, even if $h(Z\times T)$ intersects $B$.

\begin{prop}\Label{finstab}
Let $S^1$ act on $S$ with finite stabilizers.  Assume
that $A$, $A'$ and $B$ are $S^1$-invariant. Finally, assume that there
exists an $S^1$-equivariant homotopy $h:T\times Z\to S$ from $A$ to
$A'$.  Then $L_S(A',B)=L_S(A,B)$.
\end{prop}
\begin{pf}
The condition that $h$ be equivariant means that $S^1$ acts on $Z$
with finite stabilizers and 
that $h$ is equivariant, i.e. $h(t,\gamma\cdot z)=\gamma\cdot h(t,z)$,
for all $\gamma\in S^1$ and $(t,z)\in T\times Z$. 

We form the quotient orbifold $X=[S/S^1]$, which is compact of
dimension $2p$. It comes together with a
principal $S^1$-bundle $\pi:S\to X$. Let $\widetilde A$, $\widetilde
A'$, $\widetilde B$ and $\widetilde Z$ be the quotient orbifolds
obtained from $A$, $A'$, $B$ and $Z$. The homotopy $h$ descends to a
differentiable homotopy $h:T\times \widetilde Z\to X$ between
$\widetilde A$ and $\widetilde A'$, proving that $\cl(\widetilde
A)=\cl(\widetilde A')\in H^{p+1}(X)$. This conclusion is all we need
the homotopy $h$ for. 

Next we will construct, for a fixed $B$,  an element $\eta\in
H^{p-1}(X)$, such that  
$$L_S(A,B)=\int_X\eta\cup\cl(\tilde{A})\,,$$ for any $A$, such that
$A\cap B=\varnothing$. This will conclude the proof of the
proposition.

In fact, let $\beta\in H^p(S\setminus B)$, such that $\del \beta=1\in H^0(B)$. 
The $S^1$-bundle $S\setminus B\to X\setminus \widetilde B$ induces a
homomorphism 
$\pi_!:H^p(S\setminus B)\to H^{p-1}(X\setminus \widetilde B)$.  Note that the
restriction $H^{p-1}(X)\to H^{p-1}(X\setminus \widetilde B)$ is an
isomorphism, since the codimension of $\widetilde B$ in $X$ is
$p+1$. Thus, there exists a unique $\eta\in H^{p-1}(X)$, such that  
$$\eta|_{X\setminus \widetilde B}=\pi_!\beta\,.$$
Hence
$$L_S(A,B)=\int_A\beta=\int_{\widetilde A}\pi_!\beta
=\int_{\widetilde A}\eta=\int_X\eta\cup\cl(\widetilde A)\,,$$
as claimed. The last equality follows from naturality of the wrong way
maps and the projection formula.
\end{pf}

\subsection{The proof of $\nu_X(P)=(-1)^n$}

We return to the convention that dimensions are complex dimensions.

Let $X$ be a scheme with a $\Gm$-action. Let $P\in X$ be a fixed point
of this action.  The point $P$ is called an {\em   isolated }fixed
point, if 0 is not a weight of the induced action  of $\Gm$ on the
Zariski tangent space $T_X|_P$.

\begin{prop}\Label{hardwork} 
Let $M$ be a smooth scheme on which $\Gm$ acts with an isolated fixed
point $P\in M$. Let $\omega$ be an invariant (homogeneous of degree
zero) almost closed 1-form on $M$ and $X=Z(\omega)$. Assume $P\in X$.
Then
$$\nu_X(P)=(-1)^{\dim M}.$$
\end{prop}
\begin{pf}
We will use the expression of $\nu_X(P)$ as a linking number from
Proposition~4.22 of \cite{Beh}. We choose \'etale homogeneous coordinates
$x_1,\ldots,x_n$ for $M$ around $P$ and the induced \'etale coordinates
$x_1,\ldots,x_n,p_1,\ldots,p_n$  of $\Omega_M$. Since the linking
number in question is defined inside a sufficiently small sphere in
$\Omega_M$ around $P$ (and is a topological invariant), we may as well
assume that $M=\cc^n$ and $P$ is the origin. Of course, $\omega$ is
then a 1-form holomorphic (instead of algebraic) at the origin.  We
write $\omega=\sum_{i=1}^n f_i dx_i$. 

As in [ibid.], for $\eta\in\cc$, $\eta\not=0$, we write $\Gamma_\eta$
for the graph of the section 
$\frac{1}{\eta}\omega$ of $\Omega_M$.  It is defined as a subspace of
$\Omega_M$ by the equations $\eta p_i=f_i(x)$. It is oriented so that
$M\to \Gamma_\eta$ is orientation preserving.

For $t\in \rr$, we write $\Delta_t$ for the subspace of
$\Omega_M$ defined by the equations $t p_i=\ol x_i$. We orient
$\Delta_1$ in such a way that the map $\cc^n\to \Delta_1$ given by
$(x_1,\ldots,x_n)\mapsto 
(x_1,\ldots,x_n,\ol x_1,\ldots,\ol x_n)$ preserves orientation. Then
we orient all other $\Delta_t$ by continuity.  This
amounts to the same as saying that the map $(p_1,\ldots,p_n)\mapsto
(t\ol p_1,\ldots,t\ol p_n, p_1,\ldots,p_n)$ from $\cc^n$ to $\Delta_t$
preserves orientation up to a factor of $(-1)^n$. 

Proposition 4.22 of [ibid.] says that for sufficiently small
$\epsilon>0$ there exists $\eta\not=0$ such that $\Gamma_\eta'=\Gamma_\eta\cap
S_\epsilon$ is a manifold disjoint from $\Delta_1'=\Delta_1\cap S_\epsilon$ and 
$$\nu_X(P)=L_{S_\epsilon}(\Delta_1',\Gamma_\eta')\,.$$ 
Here $S_\epsilon$ is the sphere of radius
$\epsilon$ centered at the origin $P$ in $\Omega_M$. It has dimension
$4n-1$. 
Let us fix $\epsilon$ and $\eta$. 

The given $\Gm=\cc^\ast$-action on $M$ induces an action on
$\Omega_M=\cc^{2n}$.  Let us denote the degree of $x_i$ by $r_i$.
Then the degrees of $p_i$ and $f_i$ are both equal to  $-r_i$. By
restricting to $S^1\subset\cc^\ast$, we get an induced $S^1$-action on
$S_\epsilon$. This action has finite stabilizers, because none of the
$r_i$ vanish, $P$ being an isolated fixed point of the
$\Gm$-action. Note that $\Gamma_\eta'$ is an
$S^1$-invariant submanifold of $S_\epsilon$.

Consider the map from $\rr\times S^{2n-1}\to S_\epsilon$ given by
$$(t,p_1,\ldots,p_n)\mapsto \frac{\epsilon}{\sqrt{1+t^2}}(t\ol
p_1,\ldots,t\ol p_n, p_1,\ldots,p_n)\,.$$
This map is an $S^1$-equivariant homotopy between the invariant submanifolds
$\Delta_0'=\Delta_0\cap S_\epsilon$ and $\Delta_1'$.

The fact that $\Delta_1'$ is disjoint from $\Gamma_\eta'$ follows from
the fact that $\omega$ is almost closed, as explained in [ibid.].  The
fact that $\Delta_0'$ is disjoint from $\Gamma_\eta'$ is trivial:
$\Delta_0$ is (up to orientation) the fiber of the vector bundle
$\Omega_M\to M$ over the 
origin and $\Gamma_\eta$ is the graph of a section. But there is no
reason (at least none apparent to the authors)
why there shouldn't exist values of $t$ other than $0$ or $1$, for
which $\Delta_t'=\Delta_1\cap S_\epsilon$ intersects $\Gamma_\eta$. 

Still, Proposition~\ref{finstab} implies that
$$L_{S_\epsilon}(\Delta_1',\Gamma_\eta')=
L_{S_\epsilon}(\Delta_0',\Gamma_\eta')\,.$$ 
Let us denote the  fiber of $\Omega_M$ over the
origin by $\ol\Delta_0$, and its intersection with $S_\epsilon$ by
$\ol\Delta_0'$. 
By the correspondence between linking numbers and intersection numbers
(see~\cite{fulton}, Example~19.2.4), we see that
$L_{S_\epsilon}(\ol\Delta_0',\Gamma_\eta')$ is equal to the
intersection number of $\ol\Delta_0$ with $\Gamma_\eta$ at the
origin. This number is 1, as the section $\Gamma_\eta$ intersects the
fiber $\ol\Delta_0$ transversally.

Since the orientations of $\Delta_0$ and $\ol\Delta_0$ differ by
$(-1)^n$, we conclude that
$$\nu_X(P)=L_{S_\epsilon}(\Delta_1',\Gamma_\eta')=
L_{S_\epsilon}(\Delta_0',\Gamma_\eta')=
(-1)^n L_{S_\epsilon}(\ol\Delta_0',\Gamma_\eta')=(-1)^n\,,$$  
which is what we set out to prove.
\end{pf}

\begin{thm}
Let $X$ be an affine $\Gm$-scheme with an isolated fixed point $P$.
Assume that $X$ admits an equivariant symmetric obstruction theory. 
Then $$\nu_X(P)=(-1)^{\dim T_X|_P}\,.$$
\end{thm}
\begin{pf}
Let $n=\dim T_X|_P$. By Lemma~\ref{annoying}, we can assume that $X$
is embedded equivariantly in a smooth scheme $M$ of dimension $n$ and
that $X$ is the zero locus of an invariant almost closed 1-form on
$M$. Note that the embedding $X\hookrightarrow M$ identifies $T_X|_P$
with $T_M|_P$, so that $P$ is an isolated point of the $\Gm$-action on
$M$. Thus Proposition~\ref{hardwork} implies that $\nu_X(P)=(-1)^n$. 
\end{pf}

\begin{cor}\Label{weird}
Let $X$ be a $\Gm$-scheme such that all fixed points are isolated and
every fixed point admits an invariant affine open neighborhood over
which there exists an equivariant symmetric obstruction theory.  Then
we have 
$$\tilde\chi(X)=\sum_P (-1)^{\dim T_X|_P}\,,$$
the sum extending over the fixed points. Moreover, if $Z\subset X$ is
an invariant locally closed subscheme, we have
$$\tilde\chi(Z,X)=\sum_{P\in Z}(-1)^{\dim T_X|_P}\,,$$
the sum extending over the fixed points in $Z$.
\end{cor}
\begin{pf}
The product property of $\nu$ implies that $\nu_X$ is constant on
non-trivial $\Gm$-orbits.  The Euler characteristic of a scheme
on which $\Gm$ acts without fixed points is zero. These two facts
imply that only the fixed points contribute to
$\tilde\chi(X)=\chi(X,\nu_X)$.  
\end{pf}

\begin{cor}
Let $X$ be a projective scheme with a linear $\Gm$-action.
Let $X$ be endowed with an equivariant
symmetric 
obstruction theory. Assume all fixed points of $\Gm$ on $X$ are
isolated. Then we have 
$$\#^\vir(X)=\sum_P (-1)^{\dim T_X|_P}\,,$$
the sum extending over the fixed points of $\Gm$ on $X$.
\end{cor}
\begin{pf}
We use the fact that $X$ can be equivariantly embedded into a smooth
scheme to prove that every fixed point has an invariant affine open
neighborhood. Thus Corollary~\ref{weird} applies. The main result of
\cite{Beh}, Theorem~4.18, says that
$\#^\vir=\tilde\chi(X)$.
\end{pf}

\unsure{
\begin{rem} In the assumptions of the Corollary, the Theorem shows that
  for any $T$-equivariant embedding of $X$ in a smooth $M$ such that
  $x$ is isolated in $M^{T}$, the dimension of $M$ must have the same
  parity as $\dim T_xX$. 
Note that you can always increase $\dim M$ by $2$ by adding variables
$y_1,y_2$ with opposite degrees and replacing $\omega$ by
$\omega+d\,(y_1y_2)$. 
\end{rem}}

\subsubsection{Application to Lagrangian intersections}

Let $M$ be an algebraic symplectic manifold with a Hamiltonian
$\Gm$-action. Assume all fixed points are isolated. Let $V$ and $W$ be
invariant Lagrangian submanifolds, $X$ their intersection.

\begin{prop}
We have 
$$\tilde\chi(X)=\sum_{P\in X}(-1)^{\dim T_X|_P}\,,$$
the sum extending over all fixed points inside $X$.
\end{prop}
\begin{pf}
One checks that the action of $\Gm$ being Hamiltonian, i.e., that
$\Gm$ preserves the symplectic form, implies that the symmetric
obstruction theory on $X$ is equivariant symmetric. 
\end{pf}

\begin{prop}
Assume $X$ is compact. Then 
$$\deg([V]\cap[W])=\sum_{P\in X}(-1)^{\dim T_X|_P}\,,$$
the sum extending over the fixed points contained in $X$.
\end{prop}
\begin{pf}
Note that, in fact, the virtual number of points of $X$ is the
intersection number of $V$ and $W$.
\end{pf}

\begin{cor}
Assume that $X$ is compact and that $\dim T_X|_P$ is even, for all
fixed points $P$. 
Then we have
$$\deg([V]\cap[W])=\chi(X)\,.$$
\end{cor}

\eject
\section{Hilbert schemes of threefolds}

\subsection{The threefold $\A^3$}

Let $T={\mathbb G}_m^3$ be the standard 3-dimensional torus with
character group $\zz^3$.  Let $T_0$ be the kernel of the character
$(1,1,1)$. Thus,
$$T_0=\{(t_1,t_2,t_3)\in T\mid t_1t_2t_3=1\}\,.$$

We let $T$ act in the natural way on $\A^3$. Write coordinates on
$\A^3$ as $x,y,z$, then, as elements of the affine coordinate ring
$\cc[x,y,z]$ of $\A^3$, the weight of $x$ is $(1,0,0)$, the weight of
$y$ is $(0,1,0)$ and the weight of $z$ is $(0,0,1)$. 

We choose on $\A^3$ the standard 3-form $dx\wedge dy\wedge dz$ to fix
a Calabi-Yau structure.  The torus $T_0$ acts by
automorphisms of $\A^3$ preserving the Calabi-Yau structure. 

by Proposition~\ref{edt} we obtain a $T_0$-equivariant symmetric
obstruction theory on $X=\Hilb^n\A^3$.

\begin{lem}\label{keylem} 
{\rm (a)} The $T_0$-action on $X$ has a finite number of fixed
points.  These correspond to monomial ideals in $\cc[x,y,z]$.  

\noindent{\rm (b)} If $I$ is such an ideal, the $T_0$-action on the
Zariski tangent space to $X$ at $I$ has no invariant subspace.

\noindent\rm{(c)} If $I$ is such an ideal and $d$ is the
dimension of the Zariski tangent space to $X$ at $I$, we have
$(-1)^d=(-1)^n$, in other words, the integer $d$ has the same parity
as $n$.
\end{lem}
\begin{proof} 
(a) Since the $T_0$-action on $\A^3$ has the origin as unique fixed
point, any invariant subscheme must be supported at the origin. Let
$I\subset \cpx[x,y,z]$ be the corresponding ideal; $I$ must be
generated by eigenvectors of the torus action on the polynomial
ring. Any eigenvector can be written uniquely in the form $m\,g(xyz)$
where $m$ is a monomial and $g\in \cpx[t]$ is a polynomial with
$g(0)\not=0$. However, since the ideal is supported at the
origin, the zero locus of $g(xyz)$ is disjoint from the the zero locus
of $I$, and so by Hilbert's Nullstellensatz, the monomial $m$ is also
in $I$. Hence every $T_0$-invariant  ideal is generated by monomials.

(b) Let us write $A=\cc[x,y,z]$. The tangent space in question is 
$\Hom(I,A/I)$. We will prove that none of the weights
$w=(w_1,w_2,w_3)$ of $T$ on $\Hom_A(I,A/I)$ can satisfy
$w_1,w_2,w_3<0$ or $w_1,w_2,w_3\geq0$. In particular, none of these
weights can be an integer multiple of $(1,1,1)$.

This will suffice, in view of the following elementary fact:
Let $w_1,\ldots,w_n\in \zz^3$ be characters of $T$.  If none of the $w_i$
is an integer multiple of $(1,1,1)$, there exists a one-parameter
subgroup $\lambda:\Gm\hookrightarrow T_0$, such that
$w_i\circ\lambda\not=0$, for all $i=1,\ldots,n$. 

Suppose, then, that $\phi:I\to A/I$ is an eigenvector of $T$ with weight
$(w_1,w_2,w_3)$, with $w_1\geq0$, $w_2\geq0$ and $w_3\geq0$. 
Then for a monomial $x^ay^bz^c\in I$ we have
$\phi(x^ay^bz^c)\equiv x^{a+w_1}y^{b+w_2}z^{c+w_3}\mod I$, which
vanishes in $A/I$, proving that $\phi=0$. 

Now suppose $\phi:I\to A/I$ is an eigenfunction whose weights satisfy
$w_1<0$, $w_2<0$ and $w_3<0$. Let $a$ be the smallest integer such
that $x^a\in I$.  Then let $b$ be the smallest integer such that
$x^{a-1}y^b\in I$. Finally, let $c$ be the smallest integer such that
$x^{a-1}y^{b-1}z^c\in I$. Then if a monomial $x^ry^sz^t$ is in $I$,
it follows that $r\geq a$, $s\geq b$ or $t\geq c$.

We have 
$$\phi(x^ay^bz^c)=xz^c\phi(x^{a-1}y^b)\equiv
xz^cx^{a-1+w_1}y^{b+w_2}\equiv x^{a+w_1}y^{b+w_2}z^c\mod I\,.$$
We also have 
$$\phi(x^ay^bz^c)=xy\phi(x^{a-1}y^{b-1}z^c)\equiv
x^{a+w_1}y^{b+w_2}z^{c+w_3}\mod I\,.$$
We conclude that
$$x^{a+w_1}y^{b+w_2}z^c-x^{a+w_1}y^{b+w_2}z^{c+w_3}\in I\,.$$
Since the
ideal $I$ is monomial, each of these two monomials is in $I$. But the
latter one cannot be in $I$.

\unsure{
(b) Let $I\subset R:=\cpx[x,y,z]$
be the ideal of a fixed zero-dimensional subscheme; i.e., $I$ is
generated by monomials. Let $d$ be the smallest integer such that
$f_1:=x^d\in I$. Let $f_2:=x^ay^bz^c$ be any other monomial generator,
where we may assume that $a,b,c$ are minimal, and that $a<d$. Assume
we can deform $I$ to $I_\eps$ flat over $\cpx[\eps]/\eps^2$ which is
also $T$-invariant. This means that we have chosen a $T$-equivariant
map $\phi:I\to R/I$ and choose as generators of $I_\eps$ elements of
the form $f+\eps g$ where $f$ is a generator of $I$ and $g$ is a
lifting to $R$ of $\phi(f)$.\\ By Artin's criterion of
flatness\footnote{I hope this is right, I can't find my copy of Artin,
Lectures on Deformations of Singularities.} any relation between
generators in $I$ must be liftable to a relation between generators in
$I_\eps$. Consider the relation $y^bz^cf_1=x^{d-a}f_2$. By
$T$-invariance, the only possible lifting of $f_1$ is $f_1$ itself,
while $f_2$ can be lifted to $f_2 (1+\eps h((xyz)^{-1})$ where $h\in
\cpx[t]$ is any polynomial of degree $\le\min\{a,b,c\}$. Lifting the
relation would require that $f_2 h((xyz)^{-1}))$ be in the ideal
generated by $f_1$ and $f_2$, which is clearly impossible unless $h$
is a constant.}

(c) This is an immediate consequence of \cite{MNOP},
Theorem~2 in \S~4.10. In fact, this theorem states that if
$w_1,\ldots,w_d$ are the weights of $T$ on the tangent space $V$, 
$$\frac{\prod_{i=1}^d(-w_i)}{\prod_{i=1}^dw_i}=(-1)^n\,$$
inside the field of rational functions on $T$. 
\end{proof}

\begin{prop}\Label{foran}
For any $T_0$-invariant locally closed subset $Z$ of $\Hilb^n\A^3$ we have
$$\tilde\chi(Z,\Hilb^n\A^3)=(-1)^n\chi(Z)\,.$$
\end{prop}
\begin{pf}
Since there are only finitely many fixed points of $T_0$ on $X$, we
can use the fact mentioned in the proof of Lemma~\ref{keylem} to find
a one-parameter subgroup $\Gm\to 
T_0$ with respect to which all weights of all tangent spaces at all
fixed points are non-zero. Thus, all $\Gm$-fixed points are isolated.
Because $\Hilb^n\A^3$ admits an equivariant embedding into projective
space (see the proof of Proposition~\ref{edt}), every fixed point has
an invariant affine open neighborhood. 

The symmetric obstruction theory on $\Hilb^n_{(n)}\A^3$ is equivariant
symmetric with respect to the induced $\Gm$-action. We can
therefore apply Corollary~\ref{weird}.  We obtain:
$$\tilde\chi(Z,\Hilb^n\A^3)=\sum_{P\in Z} (-1)^n\,,$$
where the sum extends over fixed points $P$ contained in $Z$. Since we
also have that $\chi(Z)=\#\{P\in Z, \text{$P$ fixed}\}$, the result
follows. 
\end{pf}

Let $F_n$ denote the closed subset of $\Hilb^n\A^3$ consisting of
subschemes supported at the origin. Let $\nu_n$ be the restriction of
the canonical constructible function $\nu_{\Hilb^n\A^3}$ to
$F_n$. Thus $\tilde\chi(F_n,\Hilb^n\A^3)=\chi(F_n,\nu_n)$. Note that
all $T_0$-fixed points of $\Hilb^n\A^3$ are contained in $F_n$. 

Let $M(t)=\prod_{n=1}^\infty (1-t^n)^{-n}$ be the McMahon function. It
is the generating series for 3-dimensional partitions.  Hence, if we
write $M(t)=\sum_{n=0}^\infty p_n t^n$, then $p_n$ denotes the number
of monomial ideals $I$ in $A=\cc[x,y,z]$, such that $\dim_\cc
A/I=n$. The number $p_n$ is the number of $T_0$-fixed points in $F_n$
or $\Hilb^n\A^3$. Thus, $p_n=\chi(F_n)=\chi(\Hilb^n\A^3)$. 

\begin{cor}\Label{fnu}
We have
$$\chi(F_n,\nu_n)=(-1)^n\chi(F_n)=(-1)^n p_n\,,$$
and hence
$$\sum_{n=0}^\infty\chi(F_n,\nu_n)t^n=M(-t)\,.$$
\end{cor}

\subsection{Weighted Euler characteristics of Hilbert schemes}

Let $Y$ be a smooth threefold, and $n>0$ an integer. Consider the
Hilbert scheme of $n$ points on $Y$, denoted $\Hilb^nY$. 
The scheme $\Hilb^n Y$ is connected, smooth for $n\le
3$ and singular otherwise, and reducible for large enough $n$. 

Let us denote by $\nu_Y$ the canonical constructible function on
$\Hilb^n Y$.  Our goal is to calculate
$$\tilde\chi(\Hilb^nY)=\chi(\Hilb^nY,\nu_Y)\,.$$

Let us start with a useful general lemma on Hilbert schemes.

\begin{lem}\Label{ugl}
Let $f:Y\to Y'$ be a morphism of projective schemes and $Z\subset Y$ a
closed subscheme. Assume that $f$ is \'etale in a neighborhood of
$Z$ and that the composition $Z\to Y'$, which we will denote by $f(Z)$,
is a closed immersion of schemes. 

Let $X$ be the Hilbert scheme of $Y$ which contains $Z$ and $P$ the
point of $X$ corresponding to $Z$.  Let $X'$ be the
Hilbert scheme of $Y'$ which contains $f(Z)$. Then there exists an
open neighborhood $U$ of $P$ in $X$ and an \'etale morphism $\phi:U\to
X'$, which sends a subscheme $\tilde Z\to Y$ to the composition
$\tilde Z\to Y'$. 
\end{lem}
\begin{pf}
For the existence of the open set $U$ and the morphism $\phi$, see for
example Proposition~6.1, Chapter~I of \cite{kollar}. The fact that
$\phi$ is \'etale in a neighborhood of $P$ follows from a direct
application of the  formal criterion. 
\end{pf}

\subsubsection{The closed stratum}

We start by recalling the standard stratification of $\Hilb^nY$. The
strata are indexed by partitions of $n$.  Let
$\alpha=(\alpha_1,\ldots,\alpha_r)$ be a length $r$ partition of $n$,
i.e., $\alpha_1\ge \alpha_2\ge\ldots\ge\alpha_r\ge 1$ and
$\sum_{i=1}^r \alpha_i=n$. Let $\Hilb^n_\alpha Y$ be the locus of
subschemes whose support consists of $r$ distinct points with
multiplicities $\alpha_1,\ldots,\alpha_r$.
The closed stratum is
$\Hilb^n_{(n)} Y$.  It  corresponds to subschemes supported at a single
point.  
To fix ideas, we will endow all strata with the reduced scheme
structure.

\begin{lem}
For any threefold $Y$ there is a natural morphism 
$$\pi_Y:\Hilb^n_{(n)}Y\to Y\,.$$
\end{lem}
\begin{pf}
This is a part of the Hilbert-Chow morphism $\Hilb^nY\to S^nY$ to the
symmetric product. A proof that this is a morphism of schemes can be
found, for example, in \cite{lehn}.
\end{pf}

Note that $F_n$ is the fiber of $\pi_{\A^3}$ over the origin.

\begin{lem}\Label{translate}
We have a canonical isomorphism
\begin{equation}\Label{a3}
\Hilb^n_{(n)}\A^3=\A^3\times F_n\,.
\end{equation}
Moreover, $\nu_{\A^3}=p^\ast\nu_n$, where
$p:\Hilb^n_{(n)}\A^3\to F_n$ is the projection given by~(\ref{a3}).
\end{lem}
\begin{pf}
Consider the action of the group $\A^3$ on itself by translations. 
We get an induced action of $\A^3$ on $\Hilb^n\A^3$. 
Use this action to translate a subscheme supported at a point $P$ to a
subscheme supported at the origin. Obtain the morphism
$p:\Hilb^n_{(n)}\A^3\to F_n$ in this way. 
The product morphism $\pi_{\A^3}\times p:\Hilb^n_{(n)}\A^3\to
\A^3\times F_n$ is an isomorphism. 

It is a formal  consequence of the general properties of the canonical
constructible function, that it is constant on orbits under a group
action. This implies the claim about $\nu_{\A^3}$. 
\end{pf}

\begin{lem}\Label{UPhi}
Consider an \'etale morphism of threefolds $\phi:Y\to Y'$.

\noindent{\rm(a)} 
Let $U\subset \Hilb^n Y$ be the open subscheme parameterizing
subschemes $Z\subset Y$, which satisfy: if $P$ and $Q$ are distinct
points in the support of $Z$, then $\phi(P)\not=\phi(Q)$. 
There is an \'etale morphism $\tilde\Phi:U\to \Hilb^n Y'$ sending a
subscheme of $Y$ to its image under $\phi$. 
$$\xymatrix{
\Hilb^n_{(n)} Y\dto_{\Phi}\rto & U\dto^{\tilde\Phi} \rto & \Hilb^n Y\\
\Hilb^n_{(n)} Y'\rto &  \Hilb^n Y'}$$

\noindent{\rm(b)}
The restriction of $\tilde\Phi$ to $\Hilb^n_{(n)} Y$ induces a
cartesian diagram of schemes 
$$\xymatrix{
\Hilb^n_{(n)} Y\rto^{\Phi}\dto_{\pi_Y}\ar@{}[dr]|{\Box} &
\Hilb^n_{(n)}Y'\dto^{\pi_{Y'}}\\ Y\rto^\phi & Y'}$$
\end{lem}
\begin{pf}
The existence and \'etaleness of $\tilde\Phi$ follows immediately from
Lemma~\ref{ugl}, applied to quasi-projective covers of $Y$ and
$Y'$. Part~(b) is clear. 
\end{pf}

Let $\phi:Y\to Y'$ be an \'etale morphism with induced morphism
$\Phi:\Hilb^n_{(n)} Y\to \Hilb^n_{(n)}Y'$. By Lemma~\ref{UPhi}, the
morphism $\Phi$  extends to open neighborhoods in $\Hilb^nY$ and
$\Hilb^nY'$, respectively. The extension $\tilde\Phi$ is
\'etale. Thus, we see that 
$$\Phi^\ast(\nu_{Y'})=\nu_Y\,.$$

\begin{prop}
Every \'etale morphism $\phi:Y\to\A^3$ induces an isomorphism
$\Hilb^n_{(n)}Y=Y\times F_n$. The constructible function
$\nu_Y|_{\Hilb^n_{(n)}Y}$ is obtained by pulling back $\nu_n$ via the
induced projection $\Hilb^n_{(n)}Y\to F_n$. 
\end{prop}
\begin{pf}
Combine Lemmas \ref{translate} and \ref{UPhi}(b) with each other.
\end{pf}

\begin{cor}\Label{pt}
The morphism $\pi_Y:\Hilb^n_{(n)}Y\to Y$ is a Zariski-locally trivial
fibration with fiber $F_n$. More precisely, there exists a Zariski
open cover $\{U_i\}$ of $Y$, such that for every $i$, we have
$$(\pi_Y^{-1}(U_i),\nu_Y)=(U_i,1)\times (F_n,\nu_n)\,.$$
This is a product of schemes with constructible functions on them.
\end{cor}
\begin{pf}
Every point of $Y$ admits \'etale coordinates, defined in a Zariski
open neighborhood. 
\end{pf}

\subsubsection{Reduction to the closed stratum}

From now on the threefold $Y$ will be fixed and we denote
$\Hilb^n_\alpha Y$ by $X^n_\alpha$ and $\Hilb^nY$ by $X^n$.

\begin{lem}\Label{stratprod} 
Let $\alpha=(\alpha_1,\ldots,\alpha_r)$ be a partition of $n$. 

\noindent{\rm(a)} Let $V$ be the open subscheme of
$\prod_{i=1}^rX^{\alpha_i}$ parameterizing $r$-tuples of subschemes
with pairwise disjoint support.  Then there is a morphism
$f_\alpha:V\to X^n$ mapping $(Z_1,\ldots,Z_r)$ to $Z=\bigcup_i
Z_i$. The morphism $f_\alpha$ is \'etale. Its image $U$ is open and
contains $X^n_\alpha$. Let $Z_\alpha=f^{-1}_\alpha
X^n_\alpha$:
$$\xymatrix{
Z_\alpha\rto\dto_{\text{\rm
    Galois}}\ar@{}[dr]|{\Box} & 
V\dto^{f_\alpha}\rto & \prod_i X^{\alpha_i}\\
X^n_\alpha\rto & U\rto & X^n}$$
Moreover,  the induced morphism $Z_\alpha\to
X^n_\alpha$ is a Galois cover with Galois group $G_\alpha$, where
$G_\alpha$ is the automorphism group of the partition $\alpha$. 

\noindent{\rm (b)} The scheme $Z_\alpha$ is contained in
$\prod_iX^{\alpha_i}_{(\alpha_i)}$ and has therefore a morphism
$Z_\alpha\to Y^r$. There is a cartesian diagram
$$\xymatrix{
Z_\alpha\rto\dto\ar@{}[dr]|{\Box} &\prod_iX^{\alpha_i}_{(\alpha_i)}\dto\\
Y^r_0\rto & Y^r}$$
where $Y^{r}_0$ is the open subscheme in $Y^r$ consisting of
$r$-tuples with pairwise disjoint entries. 
\end{lem}
\begin{proof}
The existence of $f_\alpha$ and the fact that it is \'etale follows
from Lemma~\ref{ugl} applied to the \'etale map $\coprod_{i=1}^r Y\to
Y$ and the subscheme $Z_1\amalg\ldots\amalg Z_r\subset\coprod_{i=1}^r Y$. 
All other facts are also straightforward to prove. 
\end{proof}

\begin{thm}\Label{tildechi}
Let $Y$ be a smooth scheme of dimension 3. Then for all $n>0$
$$\tilde\chi(\Hilb^nY)=(-1)^n\chi(\Hilb^nY)\,.$$
This implies
$$\sum_{n=0}^\infty \tilde\chi(\Hilb^nY)t^n=M(-t)^{\chi(Y)}\,.$$
\end{thm}
\begin{pf}
By formal properties of $\tilde\chi$ as proved in \cite{Beh}, we can
calculate as follows, using Lemma~\ref{stratprod}(a):
\begin{align*}
\tilde\chi(X^n)&=\sum_{\alpha\vdash n}\tilde\chi(X^n_\alpha,X^n)\\
&=\sum_{\alpha\vdash n}\tilde\chi(X^n_\alpha,U)\\
&=\sum_{\alpha\vdash n}|G_\alpha|\,\tilde\chi(Z_\alpha,V)\\
&=\sum_{\alpha\vdash n}|G_\alpha|\,\tilde\chi\Big(Z_\alpha,
        \prod_i X^{\alpha_i}\Big)\,.
\end{align*}
By Lemma \ref{stratprod}(b), and Corollary~\ref{pt}, we have that
$Z_\alpha\to Y^{\ell(\alpha)}_0$ is  a Zariski-locally trivial fibration with
fiber $\prod_i F_{\alpha_i}$. Here we have written $\ell(\alpha)$ for
the length $r$ of the partition $(\alpha_1,\ldots,\alpha_r)$. 
We conclude:
$$\tilde\chi\Big(Z_\alpha,  \prod_i   X^{\alpha_i}\Big)
=\chi(Y^{\ell(\alpha)}_0)\,\prod_i\chi(F_{\alpha_i},\nu_{\alpha_i})$$
Together with Corollary~\ref{fnu} this gives:
\begin{equation}\Label{chh}
\tilde\chi(X^n)=(-1)^n\sum_{\alpha\vdash  n}
  |G_\alpha|\,\chi(Y^{\ell(\alpha)}_0)
  \prod_i\chi(F_{\alpha_i})\,.
\end{equation}
Using the exact same arguments with the constant function 1 in place
of $\nu$ gives the same answer, except without the sign $(-1)^n$.
This proves our first claim. The second one follows then directly from
the result of \cite{cheah}, which says that $\sum_{n=0}^\infty
\chi(\Hilb^nY)t^n=M(t)^{\chi(Y)}$.
\end{pf}

\subsection{The dimension zero MNOP conjecture}

We can now prove Conjecture 1 of \cite{MNOP}. A proof of this result
was also announced by J. Li at the workshop on Donaldson-Thomas
invariants in Urbana-Champaign in March 2005.

\begin{thm}
Let $Y$ be a projective Calabi-Yau threefold. Then, for the virtual
count of $\Hilb^nY$ with respect to the Donaldson-Thomas obstruction
theory, we have
$$\#^\vir(\Hilb^n Y)=(-1)^n\chi(\Hilb^nY)\,.$$
In other words:
$$\sum_{n=0}^\infty \#^\vir(\Hilb^n Y)\,t^n=M(-t)^{\chi(Y)}\,.$$
\end{thm}
\begin{pf}
By the main result of \cite{Beh}, Theorem~4.18, 
we have $\#^\vir(\Hilb^n Y)=\tilde\chi(\Hilb^n Y)$.  Thus the result
follows from Theorem~\ref{tildechi}.
\end{pf}

\eject


\end{document}